\documentclass[12pt]{article}
\usepackage{latexsym}
\usepackage{amsfonts}
\newtheorem{theorem}{Theorem}
\newtheorem{lemma}{Lemma}
\newtheorem{proposition}{Proposition}

\begin{document}
\def\grpc{{\mathbb C}}
\def\grpz{{\mathbb Z}}
\def\qed{{\hfill $\Box$}}
\def\vecr{{\vec r}}
\def\vecR{{\vec R}}
\title{A Survey of Results on Random Random Walks on Finite Groups}
\author{Martin Hildebrand\\Department of Mathematics and Statistics\\
University at Albany\\State University of New York\\Albany, NY 12222}
\maketitle

\begin{abstract}
A number of papers have examined
various aspects of ``random random walks'' on finite groups; the purpose of
this article is to provide a survey of this work and to 
show, bring together, and discuss some of the arguments and results in this
work. This article also provides a number of exercises. Some exercises
involve straightforward computations; others involve proving details
in proofs or extending results proved in the article. This article also 
describes some problems for further study.
\end{abstract}

\section{Introduction}

Random walks on the integers $\grpz$ are familiar to many students of
probability. (See, for example, Ch. XIV of Feller, volume 
1~\cite{fellerone}, Ch. XII of Feller, volume 2~\cite{fellertwo}, 
or Ross~\cite{ross}). Such random walks are of the form
$X_0, X_1, X_2, \dots$ where $X_0=0$ and $X_m=Z_1+\dots +Z_m$ where
$Z_1, Z_2, \dots$ are independent, identically distributed random variables
on $\grpz$. A commonly studied random walk on $\grpz$ has
$P(Z_i=1)=P(Z_i=-1)=1/2$. Various questions involving such random walks
have been well studied. For example, one may ask what is the probability
that there exists an $m>0$ such that $X_m=0$. In the example with
$P(Z_i=1)=P(Z_i=-1)=1/2$, it can be shown that this probability
is 1. (See \cite{fellerone}, p. 360.) For another example, one can use
the DeMoivre-Laplace Limit Theorem to get a good approximation of the
distribution of $X_m$ for large $m$.

One can examine random walks on sets other then $\grpz$. For instance,
there are random walks on $\grpz^2$ or $\grpz^3$. A symmetric random walk
on $\grpz^2$ has $X_0=(0,0)$ and 
$P(Z_i=(1,0))=P(Z_i=(-1,0))=P(Z_i=(0,1))=P(Z_i=(0,-1))=1/4$
while a symmetric random walk on $\grpz^3$ has $X_0=(0,0,0)$ and
$P(Z_i=(1,0,0))=P(Z_i=(-1,0,0))=P(Z_i=(0,1,0))=P(Z_i=(0,-1,0))=P(Z_i=(0,0,1))=
P(Z_i=(0,0,-1))=1/6$. It can be shown for this random walk on $\grpz^2$,
the probability that there exists an $m>0$ such that $X_m=(0,0)$ is $1$ while
for this random walk on $\grpz^3$, the probability that there exists an $m>0$
such that $X_m=(0,0,0)$ is less than $1$. (See pp. 360-361 of \cite{fellerone}
for a description and proof. Feller attributes these results to 
Polya~\cite{polya} and computation of the probability in the walk on
$\grpz^3$ to McCrea and Whipple~\cite{mccrea}.)

One can similarly look at random walks on $\grpz_n$, the integers modulo
$n$, where $n$ is a positive integer. Like the walks on the integers,
the random walk is of the form $X_0, X_1, X_2, \dots$ where $X_0=0$ and
$X_m=Z_1+\dots +Z_m$ where $Z_1, Z_2, \dots$ are i.i.d. random variables on
$\grpz_n$. One example of such a random walk has
$P(Z_i=1)=P(Z_i=-1)=1/2$ and $n$ being an odd positive integer.
This random walk on $\grpz_n$ corresponds to a finite
Markov chain which is irreducible, aperiodic, and doubly stochastic.
(For more details on this notation, see Ross~\cite{ross}.) Thus the
stationary probability for this Markov chain will be uniformly
distributed on $\grpz_n$.
(If $P$ is a probability uniformly distributed on $\grpz_n$, then
$P(a)=1/n$ for each $a\in\grpz_n$.) Furthermore, after a large enough number
of steps, the position of the random walk will be close to uniformly 
distributed on $\grpz_n$.

One may consider other probabilities on $\grpz_n$ for $Z_i$ in the random walk.
For example, on $\grpz_{1000}$, we might have $P(Z_i=0)=P(Z_i=1)=P(Z_i=10)=
P(Z_i=100)=1/4$. Again Markov chain arguments can often show that the
stationary distribution for the corresponding Markov chain will be uniformly
distributed on $\grpz_n$ and that after a large enough number of steps, the
position of the random walk will be close to uniformly distributed on
$\grpz_{1000}$. A reasonable question to ask is how large should $m$ be
to ensure that $X_m$ is close to uniformly distributed. 

One may generalize the notion of a random walk to an arbitrary finite group 
$G$. We shall suppose that the group's operation is denoted by multiplication.
The random walk will be of the form $X_0, X_1, X_2, \dots$ where $X_0$ is
the identity element of $G$, $X_m=Z_mZ_{m-1}\dots Z_2Z_1$, and $Z_1, Z_2,
\dots$ are i.i.d. random variables on $G$. (A random variable $X$ on $G$ is
such that $Pr(X=g)\ge 0$ for each $g\in G$ and $\sum_{g\in G}Pr(X=g)=1$.)
An alternate definition of a random walk on $G$ has $X_m=Z_1Z_2\dots Z_{m-1}Z_m$
instead of $X_m=Z_mZ_{m-1}\dots Z_2Z_1$. If $G$ is not abelian, then the
different definitions may correspond to different Markov chains. However,
probabilities involving $X_m$ alone do not depend on which definition we are
using.

An example of a random walk on $S_n$, the group of all permutations on
$\{1,\dots ,n\}$, has $P(Z_i=e)=1/n$ where
$e$ is the identity element and $P(Z_i=\tau)=2/n^2$ for
each transposition $\tau$. Yet again, Markov 
chain arguments can show that after a large enough number of steps, this
random walk will be close to uniformly distributed over all $n!$ permutations
in $S_n$. Again a reasonable question to ask is how large should $m$ be
to ensure that $X_m$ is close to uniformly distributed over all the
permutations in $S_n$. This problem is examined in Diaconis and
Shahshahani~\cite{diashah} and also is discussed in Diaconis~\cite{diaconis}.

A number of works discuss various random walks on finite groups. A couple of
overviews are Diaconis' monograph~\cite{diaconis} and Saloff-Coste's
survey~\cite{saloffcoste}.

In this article, we shall focus on ``random random walks'' on finite groups. To 
do so, we 
shall pick a probability for $Z_i$ at random from a set of probabilities for
$Z_i$. Then, given this probability for $Z_i$, we shall examine how close
$X_m$ is to uniformly distributed on $G$. 
To measure the distance a probability
is from the uniform distribution, we shall use the variation distance. Often,
we shall look at the average variation distance of the probability of $X_m$
from the uniform distribution; this average is over the choices for the
probabilities for $Z_i$. The next section will make these ideas more
precise.

\section{Notation}

If $P$ is a probability on $G$, we define the variation distance of $P$
from the uniform distribution $U$ on $G$ by
\[
\|P-U\|={1\over 2}\sum_{s\in G}\left|P(s)-{1\over |G|}\right|.
\]
Note that $U(s)=1/|G|$ for all $s\in G$.

{EXERCISE.} Show that $\|P-U\|\le 1$.

{EXERCISE.} Show that 
\[
\|P-U\|=\max_{A\subseteq G}|P(A)-U(A)|
\]
where $A$ ranges over all subsets of $G$. (Note that $A$ does not have to be
a subgroup of $G$.)

A random variable $X$ on $G$ is said to have probability $P$ if $P(s)=Pr(X=s)$
for each $s\in G$.

If $P$ and $Q$ are probabilities on $G$, we define the {\it convolution}
of $P$ and $Q$ by
\[
P*Q(s)=\sum_{t\in G}P(t)Q(t^{-1}s).
\]
Note that if $X$ and $Y$ are independent random variables on $G$ with
probabilities $P$ and $Q$, respectively, then $P*Q$ is the probability of the
random variable $XY$ on $G$.

If $m$ is a positive integer and $P$ is a probability on $G$, we define
\[
P^{*m}=P*P^{*(m-1)}
\]
where
\[
P^{*0}(s)=\cases{1&if $s=e$\cr 0&otherwise}
\]
with $e$ being the identity element of $G$.
Thus if $Z_1, \dots, Z_m$ are i.i.d. random variables on $G$ each with
probability $P$, then $P^{*m}$ is the probability of the random variable
$Z_mZ_{m-1}\dots Z_2Z_1$ on $G$.

For example, on $\grpz_{10}$, suppose $P(0)=P(1)=P(2)=1/3$. Then 
$P^{*2}(0)=1/9$, $P^{*2}(1)=2/9$, $P^{*2}(2)=3/9$, $P^{*2}(3)=2/9$,
$P^{*2}(4)=1/9$, and $P^{*2}(s)=0$ for the remaining elements $s\in\grpz_{10}$.
Furthermore
\begin{eqnarray*}
\|P^{*2}-U\| &=&{1\over 2}\left(\left|{1\over 9}-{1\over 10}\right|
+\left|{2\over 9}-{1\over 10}\right|+\left|{3\over 9}-{1\over 10}\right|
+\left|{2\over 9}-{1\over 10}\right|+\left|{1\over 9}-{1\over 10}\right|
\right.\\
& &\left.+\left|0-{1\over 10}\right|+\left|0-{1\over 10}\right|
+\left|0-{1\over 10}\right|+\left|0-{1\over 10}\right|
+\left|0-{1\over 10}\right|\right)\\
&=&{1\over 2}
\end{eqnarray*}

{EXERCISE.} Let $Q$ be a probability on $\grpz_{10}$. Suppose
$Q(0)=Q(1)=Q(4)=1/3$. Compute $Q^{*2}$ and $\|Q^{*2}-U\|$.

{EXERCISE.} Consider the probability $P$ in the previous example and the
probability $Q$ in the previous exercise. Compute $P^{*4}$, $Q^{*4}$,
$\|P^{*4}-U\|$, and $\|Q^{*4}-U\|$.

Let $a_1, a_2, \dots a_k\in G$. Suppose $p_1, \dots, p_k$ are positive
numbers which sum to $1$. Let
\[
P_{a_1,\dots,a_k}(s)=\sum_{b=1}^kp_b\delta_{s,a_b}
\]
where
\[
\delta_{s,a_b}=\cases{1&if $s=a_b$\cr 0&otherwise.}
\]
The random walk $X_0, X_1, X_2, \dots$ where $Z_1, Z_2, \dots$ are i.i.d.
random variables on $G$ with probability $P_{a_1,\dots,a_k}$ is said to be
supported on $(a_1,\dots,a_k)$.

If we know $k$ and $p_1,\dots,p_k$, then in theory we can find
\[
\|P_{a_1,\dots,a_k}^{*m}-U\|
\]
for any $k$-tuple $(a_1,\dots,a_k)$ and positive integer $m$. Thus if we
have some probability distribution for $(a_1,\dots,a_k)$, we can in theory find
\[
E(\|P_{a_1,\dots,a_k}^{*m}-U\|)
\]
for each positive integer $m$ since this variation distance is a function of 
the random $k$-tuple $(a_1,\dots,a_k)$.

In a random random walk, a typical probability for $Z_i$ will 
be $P_{a_1,\dots,a_k}$ where $p_1=p_2=\dots=p_k=1/k$ and
$(a_1,a_2,\dots,a_k)$ chosen uniformly over all $k$-tuples
with distinct elements of $G$. However, other probabilities for $Z_i$
sometimes may be considered instead.

For example, on $\grpz_5$, suppose we let $p_1=p_2=p_3=1/3$ and we
choose $(a_1,a_2,a_3)$ at random over all $3$-tuples with distinct elements
of $\grpz_5$. Then
\begin{eqnarray*}
E(\|P_{a_1,a_2,a_3}^{*m}-U\|)&=&{1\over 10}\left(\|P_{0,1,2}^{*m}-U\|
+\|P_{0,1,3}^{*m}-U\|+\|P_{0,1,4}^{*m}-U\|\right.
\\
& &\left. +\|P_{0,2,3}^{*m}-U\|+\|P_{0,2,4}^{*m}-U\|+\|P_{0,3,4}^{*m}-U\|
\right.
\\
& &\left. +\|P_{1,2,3}^{*m}-U\|+\|P_{1,2,4}^{*m}-U\|+\|P_{1,3,4}^{*m}-U\|
\right.
\\
& &\left. +\|P_{2,3,4}^{*m}-U\|\right).
\end{eqnarray*}
Note that we are using facts such as $P_{0,1,2}=P_{2,1,0}=P_{1,2,0}$
since $p_1=p_2=p_3=1/3$; thus we averaged over $10$ terms instead of 
over $60$ terms.

We shall often write $E(\|P^{*m}-U\|)$ instead of 
$E(\|P_{a_1,\dots,a_k}^{*m}-U\|)$.

Often times we shall seek upper bounds on $E(\|P^{*m}-U\|)$. Note that
by Markov's inequality, if $E(\|P^{*m}-U\|)\le u$, then
$Pr(\|P^{*m}-U\|\ge cu)\le 1/c$ for $c>1$ where the probability is over the
same choice of the $k$-tuple $(a_1,\dots,a_k)$ as used in determining
$E(\|P^{*m}-U\|)$.

Lower bounds tend to be found on $\|P^{*m}-U\|$ for all $a_1,\dots,a_k$, and
we turn our attention to some lower bounds in the next section.

{EXERCISE.} Compute $E(\|P^{*4}-U\|)$ if $G=\grpz_{11}$, $k=3$, 
$p_1=p_2=p_3=1/3$, and $(a_1,a_2,a_3)$ are chosen uniformly from all $3$-tuples
with distinct elements of $G$. It may be very helpful to use a
computer to deal with the very tedious calculations; however, it
may be instructive to be aware of the different values of $\|P^{*4}-U\|$
for the different choices of $a_1$, $a_2$, and $a_3$. In particular, what can
you say about $\|P_{a_1,a_2,a_3}^{*4}-U\|$ and 
$\|P_{0,a_2-a_1,a_3-a_1}^{*4}-U\|$?

Throughout this article, there are a number of references to logarithms.
Unless a base is specified, the expression $\log$ refers to the
natural logarithm. 

\section{Lower bounds}
\subsection{Lower bounds on $\grpz_n$ with $k=2$}

In examining our lower bounds, we shall first look at an elementary case.
Suppose $G=\grpz_n$. We shall consider random walks supported on $2$
points. The lower bound is given by the following:

\begin{theorem}
\label{lowerboundtwopts}
Suppose $p_1=p_2=1/2$. Let $\epsilon>0$ be given. There exists values
$c>0$ and $N>0$ such that if $n>N$ and $m<cn^2$, then
\[
\|P_{a_1,a_2}^{*m}-U\|>1-\epsilon
\]
for all $a_1$ and $a_2$ in $G$ with $a_1\ne a_2$.
\end{theorem}

{\bf Proof:} Let $F_m=|\{i: 1\le i\le m, Z_i=a_1\}|$, and
let $S_m=|\{i: 1\le i\le m, Z_i=a_2\}|$. In other words, if we 
perform $m$ steps of the random walk supported on $(a_1,a_2)$
on $\grpz_n$, we add $a_1$ exactly
$F_m$ times, and we add $a_2$ exactly $S_m$ times. Note that
$F_m+S_m=m$ and that $X_m=F_ma_1+S_ma_2$.
Observe that $E(F_m)=m/2$. Furthermore, observe that by the DeMoivre-Laplace
Limit Theorem (see Ch. VII of \cite{fellerone}, for example), there are
constants $z_1$, $z_2$, and $N_1$ such that if $n>N_1$, then
$Pr(m/2-z_1\sqrt{m/4}<F_m<m/2+z_2\sqrt{m/4})>1-\epsilon/2$.
Furthermore, there exists a constant $c$ such that if $m<cn^2$ and $n>N_1$,
then $(\epsilon/6)n\ge \sqrt{m/4}$ and
$Pr(m/2-(\epsilon/6)n<F_m<m/2+(\epsilon/6)n)>1-\epsilon/2$.
Let $A_m=\{F_ma_1+(m-F_m)a_2: m/2-(\epsilon/6)n<F_m<m/2+(\epsilon/6)n\}$.
Thus $P^{*m}(A_m)>1-(\epsilon/2)$. However, $|A_m|\le
(\epsilon/3)n+1$. Thus if $n>N_1$ and $m<cn^2$, then
\begin{eqnarray*}
\|P^{*m}-U\|&>&1-{\epsilon\over 2}-\left({\epsilon\over 3}+{1\over
n}\right)\\
&=&1-{5\epsilon\over 6}-{1\over n}
\end{eqnarray*}
Let $N=\max(N_1,6/\epsilon)$. If $n>N$ and $m<cn^2$, then $\|P^{*m}-U\|>
1-\epsilon$. \qed

{EXERCISE.} Modify the above theorem and its proof to consider the case
where $0<p_1<1$ is given and $p_2=1-p_1$.

\subsection{Lower bounds on abelian groups with $k$ fixed}

Now let's turn our attention to a somewhat more general result. This result is
a minor modification of a result of Greenhalgh~\cite{greenhalgh}. Here the
probability will be supported on elements $a_1, a_2, \dots, a_k$. One 
similarity to the
proof of the previous theorem is that we use the number of times in
the random walk we choose $Z_i$ to be $a_1$, the number of times
we choose $Z_i$ to be $a_2$, etc.

\begin{theorem}
\label{greenhalghlowerbound}
Let $p_1=p_2=\dots=p_k=1/k$. Suppose $G$ is an abelian group of order $n$.
Assume $n\ge k\ge 2$.
Let $\delta$ be a given value under $1/2$. Then there exists a value 
$\gamma>0$
such that if $m<\gamma n^{2/(k-1)}$, then $\|P_{a_1,\dots,a_k}^{*m}-U\|>
\delta$ for any $k$ distinct values $a_1, a_2, \dots a_k\in G$.
\end{theorem}

{\bf Proof:} This proof slightly modifies an argument in 
Greenhalgh~\cite{greenhalgh}. That proof, which only considered the case
where $G=\grpz_n$, had $\delta=1/4$ but considered a
broader range of $p_1, p_2, \dots, p_k$.

Since $G$ is abelian, $X_m$ is completely determined by the number
of times $r_i$ we pick each element $a_i$. Observe that
$r_1+\dots +r_k=m$. Thus if $g=r_1a_1+\dots+r_ka_k$,
we get
\[
P^{*m}(g)\ge Q_m(\vecr):={m!\over\prod_{i=1}^k(r_i!k^{r_i})}
\]
where $\vecr=(r_1,\dots,r_k)$. (Note that we may use $P^{*m}$ instead
of $P_{a_1,\dots,a_k}^{*m}$.)

Suppose $|r_i-m/k|\le\alpha_i\sqrt{m}$ for $i=1,\dots,k$. Then
\begin{eqnarray*}
Pr(|r_i-m/k|\le\alpha_i\sqrt{m}, i=1,\dots,k)
&\ge&1-\sum_{i=1}^k{m(1/k)(1-1/k)\over \alpha_i^2m}\\
&=&1-\sum_{i=1}^k{(1/k)(1-1/k)\over\alpha_i^2}\\
&>&2\delta
\end{eqnarray*}
where we assume $\alpha_i$ are large enough constants to make the last
inequality hold.

Let $\vecR=\{\vecr: r_1+\dots+r_k=m {\rm \ and\ }|r_i-m/k|\le\alpha_i\sqrt m$ for $i=1,\dots,k\}$, and
let $A_m=\{r_1a_1+\dots+r_ka_k: (r_1,\dots,r_k)\in\vecR\}$.
Thus $P^{*m}(A_m)>2\delta$.

Observe that, with $\epsilon_m(\vecr)\rightarrow 0$ uniformly over
$\vecr\in\vecR$ as $m\rightarrow\infty$, we have, by Stirling's formula,
\begin{eqnarray*}
Q_m(\vecr)m^{(k-1)/2}&=&
{m^{(k-1)/2}e^{-m}m^m\sqrt{2\pi m}(1+\epsilon_m(\vecr))\over
\prod_{i=1}^kk^{r_i}e^{-r_i}r_i^{r_i}\sqrt{2\pi r_i}}\\
&=&{e^{-m}\over \exp(-\sum_{i=1}^kr_i)}\cdot
{m^m\over \prod_{i=1}^k(kr_i)^{r_i}}\cdot{1\over \sqrt{\prod_{i=1}^k\left(
{r_i\over m}\right)}}
\cdot{1+\epsilon_m(\vecr)\over (2\pi)^{(k-1)/2}}\\
&=&1\cdot{1 \over \prod_{i=1}^k(kr_i/m)^{r_i}}\cdot
{1\over\prod_{i=1}^k\sqrt{r_i/m}}
\cdot{1+\epsilon_m(\vecr)\over (2\pi)^{(k-1)/2}}
\end{eqnarray*}
provided that all of the values $r_1, \dots, r_k$ are positive integers.

It can be shown that $\prod_{i=1}^k\sqrt{r_i/m}\rightarrow(1/k)^{k/2}$
uniformly over $\vecr\in\vecR$ as $m\rightarrow\infty$.

Now let $x_i=r_i-m/k$; thus $kr_i=kx_i+m$. Note that $x_1+\dots+x_k=0$
since $r_1+\dots+r_k=m$. Observe that
\begin{eqnarray*}
\prod_{i=1}^k\left({kr_i\over m}\right)^{r_i}&=&
\prod_{i=1}^k\left(1+{kx_i\over m}\right)^{m/k+x_i}\\
&=&\exp\left(\sum_{i=1}^k\left({m\over k}+x_i\right)\log\left(1+
{kx_i\over m}\right)\right).
\end{eqnarray*}
Now observe
\begin{eqnarray*}
\lefteqn{\sum_{i=1}^k\left({m\over k}+x_i\right)
\log\left(1+{kx_i\over m}\right)}
\hspace{0.5in}\\
&=&\sum_{i=1}^k\left({m\over k}+x_i\right)\left({kx_i\over m}
-{k^2x_i^2\over 2m^2}\right)+f(m,x_1,\dots,x_k)\\
&=&{k \over 2m}\sum_{i=1}^kx_i^2-{k^2\sum_{i=1}^kx_i^3\over 2m^2}
+f(m,x_1,\dots,x_k)
\end{eqnarray*}
where for some constant $C_1>0$, $|f(m,x_1,\dots,x_k)|\le C_1/\sqrt{m}$
for all $\vecr\in\vecR$. So for some constant $C_2>0$,
\[
\left|f(m,x_1,\dots,x_k)-{k^2\sum_{i=1}^kx_i^3\over 2m^2}\right|\le
{C_2\over \sqrt{m}}
\]
for all $\vecr\in\vecR$.

Thus for some constant $\alpha>0$ and some integer $M>0$, we have
\[
Q_m(\vecr)m^{(k-1)/2}\ge \alpha
\]
for all $\vecr\in\vecR$
and $m\ge M$. (We may also assume that $M$ is large enough that
$m/k>\alpha_i\sqrt m$ if $m\ge M$. Thus if $m\ge M$, we have $r_i>0$
for all $\vecr\in\vecR$.)
Now suppose that
\[
{\alpha\over m^{(k-1)/2}}\ge{2\over n}
\]
Thus if $g\in A_m$, $P^{*m}(g)\ge\alpha/m^{(k-1)/2}\ge 2/n$
and $P^{*m}(g)-(1/n)\ge P^{*m}(g)/2$. Thus
$P^{*m}(A_m)-U(A_m)\ge 0.5P^{*m}(A_m)>\delta$, and so
\[
\|P_{a_1,\dots,a_k}^{*m}-U\|>\delta
\]
if $m\le (\alpha/2)^{2/(k-1)}n^{2/(k-1)}$ and $m\ge M$. The following exercise
shows that the above reasoning suffices if $M\le (\alpha/2)^{2/(k-1)}
n^{2/(k-1)}$, i.e. $n\ge(2/\alpha)M^{(k-1)/2}$. 

{EXERCISE.} If $P$ and $Q$ are probabilities on $G$, show that
$\|P*Q-U\|\le\|P-U\|$. Use this result to show that if $m_2\ge m_1$, 
then $\|P^{*m_2}-U\|\le\|P^{*m_1}-U\|$.

For smaller
values of $n\ge k$, observe that $\|P^{*0}-U\|=1-1/n>\delta$, and we
may choose $\gamma_1$ such that $\gamma_1 N^{2/(k-1)}<1$ where $N=
(2/\alpha)M^{(k-1)/2}$. Let $\gamma=\min(\gamma_1,(\alpha/2)^{2/(k-1)})$.
Note that this value $\gamma$ does not depend on which abelian group
$G$ we chose.
\qed

{EXERCISE.} Extend the previous theorem and proof to deal with the
case where $p_1, \dots, p_k$ are positive numbers which sum to $1$.

{EXERCISE.} Extend the previous theorem and proof to deal with the
case where $\delta$ is a given value under $1$. (Hint: If $\delta<1/b$ where
$b>1$, replace ``$2\delta$'' by ``$b\delta$'' and replace ``$2/n$''
by an appropriate multiple of $1/n$.)

\subsection{Some lower bounds with $k$ varying with the order of the group}

Now let's look at some lower bounds when $k$ varies with $n$. The following
theorem is a rephrasing of a result of Hildebrand~\cite{mvhptrf}.
\begin{theorem}
Let $G$ be an abelian group of order $n$. If $k=\lfloor(\log n)^a\rfloor$
where $a<1$ is a constant, then for each fixed positive value $b$, there
exists a function $f(n)$ such that $f(n)\rightarrow 1$ as $n\rightarrow\infty$
and $\|P_{a_1,\dots,a_k}^{*m}-U\|\ge f(n)$ for each $(a_1,\dots,a_k)\in G^k$
where $m=\lfloor (\log n)^b\rfloor$.
\end{theorem}

{\bf Proof:} The proof also is based on material in \cite{mvhptrf}.
In the first $m$ steps of the
random walk, each value $a_i$ can be picked either $0$
times, $1$ time, \dots, or $\lfloor(\log n)^b\rfloor$ times. Since the group
is abelian, the value $X_m$ depends only on the number of times each $a_i$ is 
picked in the first $m$ steps of the random walk.
So after $m$ 
steps, there are at most $(1+\lfloor (\log n)^b\rfloor)^k$ different
possible values for $X_m$. The following exercise implies the theorem.

{EXERCISE.} Prove the following. If $k=\lfloor(\log n)^a\rfloor$ and $a<1$
is a constant, then
\[
{(1+\lfloor(\log n)^b\rfloor)^k\over n}\rightarrow 0
\]
as $n\rightarrow\infty$.

Note that the function $f$ may depend on $b$ but does not necessarily depend on
which abelian group $G$ we chose.
\qed

The following lower bound, mentioned in \cite{mvhptrf},
applies for any group $G$.
\begin{theorem}
Let $G$ be any group of order $n$. Suppose $k$ is a function of $n$.
Let $\epsilon>0$ be given. If
\[
m=\lfloor {\log n\over \log k}\rfloor(1-\epsilon),
\]
then $\|P_{a_1,\dots,a_k}^{*m}-U\|\ge f(n)$ for each $(a_1,\dots,a_k)\in G^k$
and some $f(n)$ such that $f(n)\rightarrow 1$ as $n\rightarrow\infty$.
\end{theorem}

{\bf Proof:} Note that $X_m$ has at most $k^m$ possible values and that
$k^m\le n^{1-\epsilon}$. Thus $\|P_{a_1,\dots,a_k}^{*m}-U\|\ge
1-(n^{1-\epsilon}/n)$. Let $f(n)=1-(n^{1-\epsilon}/n)$.
\qed

If $k=\lfloor (\log n)^a\rfloor$ and $a>1$ is constant, then 
the previous lower bound can be
made slightly larger for abelian groups. The following theorem is
a rephrasing of a result in Hildebrand~\cite{mvhptrf}.

\begin{theorem}
Let $G$ be an abelian group of order $n$. Suppose $k=\lfloor(\log n)^a\rfloor$
where $a>1$ is a constant. Let $\epsilon>0$ be given. Suppose
$p_1=\dots=p_k=1/k$. Suppose
\[
m=\lfloor{a\over a-1}{\log n\over\log k}(1-\epsilon)\rfloor.
\]
Then for some $f(n)$ such that $f(n)\rightarrow 1$ as $n\rightarrow\infty$, 
$\|P_{a_1,\dots,a_k}^{*m}-U\|\ge f(n)$ for all $(a_1,\dots,a_k)$ such that
$a_1, \dots, a_k$ are distinct elements of $G$.
\end{theorem}

{\bf Proof:} The proof of this theorem is somewhat trickier than the
previous couple of proofs and is based on a proof in \cite{mvhptrf}. 
We shall find functions $g(n)$ and $h(n)$
such that $g(n)\rightarrow 0$ and $h(n)\rightarrow 0$ as $n\rightarrow
\infty$ and the following holds. Given $(a_1,\dots,a_k)\in G^k$, there
exists a set $A_m$ such that $P_{a_1,\dots,a_k}^{*m}(A_m)>1-g(n)$ while
$U(A_m)<h(n)$.

To find such a set $A_m$, we use the following proposition.
\begin{proposition}
Let $R=\{j: 1\le j\le m$ and $Z_j=Z_i$ for some $i<j\}$. Then there exist
functions $f_1(n)$ and $f_2(n)$ such that $f_1(n)\rightarrow 0$ as 
$n\rightarrow\infty$, $f_2(n)\rightarrow 0$ as $n\rightarrow \infty$,
and $Pr(|R|>f_1(n)m)<f_2(n)$.
\end{proposition}

{\bf Proof:} First note $Pr(Z_j=Z_i$ for some $i<j)\le (j-1)/k$. Thus
$\sum_{j=1}^m Pr(Z_j=Z_i$ for some $i<j)<m^2/k$. Thus $E(|R|)<m^2/k$.
Thus by Markov's inequality,
\[
Pr(|R|>f_1(n)m)<{m^2/k\over f_1(n)m}={m\over f_1(n)k}
\]
for any function $f_1(n)>0$. Since $m/k\rightarrow 0$ as $n\rightarrow
\infty$, we can find a function $f_1(n)\rightarrow 0$ as $n\rightarrow \infty$
such that $f_2(n):=m/(f_1(n)k)\rightarrow 0$ as $n\rightarrow\infty$. \qed

Let $A_m=\{Z_mZ_{m-1}\dots Z_1$ such that $|\{j: 1\le j \le m$ and
$Z_j=Z_i$ for some $i<j\}|\le f_1(n)m\}$. So by construction,
\[
P_{a_1,\dots,a_k}^{*m}(A_m)\ge 1-f_2(n).
\]

Now let's consider $U(A_m)$. There are $k^m$ different choices overall
for $Z_1, Z_2, \dots, Z_m$. Observe that for each $X_m\in A_m$, there is at 
least one way to obtain it after $m$ steps of the random walk such 
that the values $Z_1, \dots, Z_k$
have $m-\lceil f_1(n)m\rceil$ distinct values. Rearranging these
distinct values does not change the value $X_m$ since $G$ is abelian.
Thus each $X_m$ is obtained by at least $(m-\lceil f_1(n)m\rceil)!$
different choices in the walk. Thus $A_m$ has at most 
$k^m/(m-\lceil f_1(n)m\rceil)!$ different values. The following 
exercise completes the proof of the theorem.

{EXERCISE.} Show that for the values $m$ and $k$ in this theorem
and for any function $f_1(n)>0$ such that $f_1(n)\rightarrow 0$ as
$n\rightarrow \infty$,
\[
{k^m/(m-\lceil f_1(n)m\rceil)! \over n}\rightarrow 0
\]
as $n\rightarrow \infty$.
(Hint: Use Stirling's formula to show that $(m-\lceil f_1(n)m\rceil)!$
is $n^{(1-\epsilon)(1-g_1(n))/(a-1)}$ for some function $g_1(n)\rightarrow 0$
as $n\rightarrow\infty$, and search for terms which are
insignificant when one expresses $k^m$ in terms of $n$.)

\qed

Note that the function $f(n)$ in the previous theorem does not have to depend
on the which group $G$ of order $n$ we have.

\section{Upper bounds}

We now turn our attention to upper bounds for the random random walks on 
finite groups. Usually these bounds will refer to the expected value of
the variation distance after $m$ steps. This expected value, as noted in
an earlier section, is over the choice of the $k$-tuple $(a_1,\dots,a_k)$.
First we shall look at a lemma frequently used for upper bounds involving
random walks on finite groups as well as random random walks on finite groups.

\subsection{The Upper Bound Lemma of Diaconis and Shahshahani}

This lemma uses techniques from Fourier analysis on finite groups.
A more extensive summary of this technique appears in 
Diaconis~\cite{diaconis}. There are a number of sources which
describe Fourier analysis on finite groups (e.g. Terras~\cite{terras})
and representation theory of finite groups (e.g. Serre~\cite{serre}
or Simon~\cite{simon}).

A representation $\rho$ of a finite group $G$ is a function from $G$
to $GL_n(\grpc)$ such that $\rho(st)=\rho(s)\rho(t)$ for all $s, t \in G$;
the value $n$ is called the degree of the representation and is denoted 
$d_{\rho}$. For example, if $j\in\{0, 1, \dots, n-1\}$, then
$\rho_j(k)=\left[e^{2\pi ijk/n}\right]$ for $k\in \grpz_n$ is a
representation of $\grpz_n$. For any group $G$, the representation
$\rho(s)=[1]$ for all $s\in G$ is called the trivial representation.
A representation is said to be
irreducible if there is no proper nontrivial subspace $W$ of $\grpc^n$
(where $n=d_{\rho}$) such that $\rho(s)W\subseteq W$ for all $s\in G$.
If there exists an invertible complex matrix $A$ such that $A\rho_1(s)A^{-1}
=\rho_2(s)$ for all $s\in G$, then the representation $\rho_1$ and
$\rho_2$ are said to be equivalent. It can be shown that each irreducible
representation is equivalent to a unitary representation. We shall assume
that when we pick an irreducible representation up to equivalence, we pick
a unitary representation.
It can be shown that $|G|=\sum_{\rho}d_{\rho}^2$ where the sum is over
all irreducible representations of $G$ up to equivalence.

We define the Fourier transform 
\[
\hat P(\rho)=\sum_{s\in G}P(s)\rho(s).
\]
The following lemma, known as the Upper Bound Lemma, is due to Diaconis
and Shahshahani~\cite{diashah}
 and is frequently used in studying probability on finite 
groups. The description here is based on the description in 
Diaconis~\cite{diaconis}.
\begin{lemma}
Let $P$ be a probability on a finite group $G$ and $U$ be the uniform
distribution on $G$. Then
\[
\|P-U\|^2\le{1\over 4}\sum_{\rho}^{*}d_{\rho}Tr(\hat P(\rho)\hat
P(\rho)^{*})
\]
where the sum is over all non-trivial irreducible representations $\rho$ up
to equivalence and $*$ of a matrix denotes its conjugate transpose.
\end{lemma}

\subsection{Upper bounds for random random walks on $\grpz_n$ where $n$ is
prime}

A result shown in Hildebrand~\cite{mvhptrf} (and based upon a result in
\cite{mvhthesis}) is the following.
\begin{theorem}
\label{mvhptrfthm}
Suppose $k$ is a fixed integer which is at least $2$. Let $p_i, i=1,\dots,k$
be such that $p_i>0$ and $\sum_{i=1}^kp_i=1$. Let $\epsilon>0$ be given.
Then for some values $N$ and $\gamma>0$ (where $N$ and $\gamma$ may depend
on $\epsilon, k, p_1, \dots, p_k$, but not $n$), 
\[
E(\|P_{a_1,\dots,a_k}^{*m}-U\|)<\epsilon 
\]
for $m=\lfloor\gamma n^{2/(k-1)} \rfloor$ 
for prime numbers $n>N$. The expectation is over a uniform choice
of $k$-tuples $(a_1, \dots, a_k)\in G^k$ such that $a_1, \dots, a_k$
are all distinct.
\end{theorem}

{\bf Proof:} This presentation is based upon the ideas in the proof in
\cite{mvhptrf}.
First the result of the following exercise means that we may use
the Upper Bound Lemma.

{EXERCISE.} Suppose that given $\epsilon^{\prime}>0$, there exists values 
$\gamma^{\prime}>0$ and $N^{\prime}$ (which may depend 
on $\epsilon^{\prime}, k, p_1, \dots p_k$, but 
not on $n$) such that 
\[
E(\|P_{a_1, \dots, a_k}^{*m}-U\|^2)<\epsilon^{\prime}
\]
if $m=\lfloor\gamma^{\prime} n^{2/(k-1)}\rfloor$ and 
$n$ is a prime which is greater
than $N^{\prime}$ where the expectation is as in Theorem~\ref{mvhptrfthm}. Then 
Theorem~\ref{mvhptrfthm} holds.

The following proposition is straightforward, and its proof is left to the
reader. Note that by abuse of notation, we view the Fourier transform in
this proposition as a scalar instead of as a $1$ by $1$ matrix.

\begin{proposition}
\[
|\hat P_{a_1,\dots,a_k}(j)|^2=\left(\sum_{i=1}^kp_i^2\right)+
2\sum_{1\le i_1<i_2\le k}p_{i_1}p_{i_2}\cos(2\pi(a_{i_1}-a_{i_2})j/n)
\]
where $\hat P_{a_1,\dots,a_k}(j)=\hat P_{a_1,\dots,a_k}(\rho_j)$.
\end{proposition}

In the rest of the proof of Theorem~\ref{mvhptrfthm}, we shall assume $j\ne 0$
since  the $j=0$ term corresponds to the trivial representation, which is
not included in the sum in the Upper Bound Lemma.

Let's deal with the case $k=2$ now. We see that
\[
\hat P_{a_1,a_2}(j)=(1-2p_1p_2)+2p_1p_2\cos(2\pi(a_1-a_2)j/n).
\]
Note that $(a_1-a_2)j$ mod $n$ runs through $1, 2, \dots, n-1$.
Thus
\begin{eqnarray*}
\sum_{j=1}^{n-1}|\hat P_{a_1,a_2}(j)|^{2m}&=&
\sum_{j=1}^{n-1}(1-2p_1p_2+2p_1p_2\cos(2\pi j/n))^m\\
&\le &2\sum_{j=1}^{\lfloor(n-1)/2\rfloor}\exp(-cj^2m/n^2)\\
&\le&2\sum_{j=1}^{\infty}\exp(-cjm/n^2)\\
&=&2{\exp(-cm/n^2)\over 1-\exp(-cm/n^2)}
\end{eqnarray*}
for some constant $c>0$. (This argument is similar to one in
Chung, Diaconis, and Graham~\cite{cdg}.) For some
$\gamma>0$, if $m=\lfloor \gamma n^2\rfloor$, then
$\|P_{a_1,a_2}^{*m}-U\|<\epsilon$ for sufficiently large
primes $n$ uniformly over all $a_1, a_2\in\grpz_n$ with $a_1\ne a_2$.

From now on in the proof of Theorem~\ref{mvhptrfthm}, we assume $k\ge 3$.
Note that $|\hat P_{a_1,\dots,a_k}(j)|^2\ge 0$ and that if
$\cos(2\pi(a_{i_1}-a_{i_2})j/n)\le 0.99$ for some $i_1$ and $i_2$ with
$1\le i_1<i_2\le k$, then $|\hat P_{a_1,\dots,a_k}(j)|^2\le b_1:=
1-0.02\min_{i_1\ne i_2}p_{i_1}p_{i_2}$.

{EXERCISE.} If $m=\lfloor \gamma n^{2/(k-1)}\rfloor$ where $\gamma>0$
is a constant, show that $\lim_{n\rightarrow \infty}nb_1^m=0$.

Thus we need to focus on values of $a_1,\dots,a_k$ such that
$\cos(2\pi(a_{i_1}-a_{i_2})j/n)>0.99$ for all $i_1$ and $i_2$ with
$1\le i_1<i_2\le k$. In doing so, we'll focus on the case where
$i_1=1$.

Let $g_n(x)=x_0$ where $x_0\in(-n/2,n/2]$ and $x\equiv x_0\pmod n$.
Thus $\cos(2\pi x/n)=\cos(2\pi g_n(x)/n)$.
The following lemma looks at the probability that
\[
(g_n((a_1-a_2)j)/n,\dots,g_n((a_1-a_k)j)/n)
\]
falls in a given $(k-1)$-dimensional ``cube'' of a given size.

\begin{lemma}
\label{cubelemma} Suppose $k\ge 3$ is constant.
Let $\epsilon>0$ be given. Then there exists a value $N_0$ such that
if $n>N_0$ and $n$ is prime, then
\begin{eqnarray*}
\lefteqn{P\left(m_i(\epsilon/2)^{1/(k-1)}n^{(k-2)/(k-1)}/(2n) \le
g_n((a_1-a_i)j)/n \right.}\hspace{0.5in}\\
&&\left.
\le (m_i+1)(\epsilon/2)^{1/(k-1)}n^{(k-2)/(k-1)}/(2n), i=2,\dots k\right)
\\
&\le&{1.1(\epsilon/2)\over 2^{k-1}n}
\end{eqnarray*}
for each $(m_2,\dots,m_k)\in\grpz^{k-1}$ and $j\in\{1,\dots,n-1\}$.
\end{lemma}

{\bf Proof:} Observe that if $n$ is prime and odd, then $g_n((a_1-a_i)j)$ 
may be any of the values $(-n+1)/2,\dots,-2,-1,1,2,\dots,(n-1)/2$ except
for those values taken by $g_n((a_1-a_{i^{\prime}})j)$ for $i^{\prime}<i$.
Furthermore, different values of $a_i$ correspond to different values of
$g_n((a_1-a_i)j)$. Note that this statement need not hold if $n$ were not prime
or if $j$ were $0$.

On the interval $[a,b]$, there are at most $b-a+1$ integers. Thus on the 
interval
\[
\left[{m_i(\epsilon/2)^{1/(k-1)}n^{(k-2)/(k-1)}\over 2},
{(m_i+1)(\epsilon/2)^{1/(k-1)}n^{(k-2)/(k-1)}\over 2}\right],
\]
there are at most $1+(\epsilon/2)^{1/(k-1)}n^{(k-2)/(k-1)}/2$ 
possible values of $g_n((a_1-a_i)j)$.

Thus the probability in the statement of the lemma is less than or
equal to
\[
{\left(1+(\epsilon/2)^{1/(k-1)}n^{(k-2)/(k-1)}/2\right)^{k-1}
\over (n-1)(n-2)\dots (n-k+1)}\sim {(\epsilon/2)\over 2^{k-1}n}
\]
where $f(n)\sim g(n)$ means $\lim_{n\rightarrow \infty} f(n)/g(n)=1$.
The lemma follows. \qed

{EXERCISE.} Explain why the proof of the previous lemma fails if $k=2$
(even though this failure is not acknowledged in \cite{mvhptrf}).

Next we wish to find an upper bound on $|\hat P_{a_1,\dots,a_k}(j)|^{2m}$
for $(g_n((a_1-a_2)j)/n,\dots,g((a_1-a_k)j)/n)$ in each such 
$(k-1)$-dimensional ``cube''. Note that there is a constant $c_1\in (0,1]$
such that if $\cos(2\pi k_1/n)>0.99$, then $\cos(2\pi k_1/n)\le 
1-(c_1/2)(g_n(k_1))^2/n^2$. Now let $\ell_n$ be largest positive integer $\ell$
such that
\[
\cos\left(2\pi\left({\ell(\epsilon/2)^{1/(k-1)}n^{(k-2)/(k-1)}\over 2n}\right)
\right)>0.99
\]
and 
\[
{\ell(\epsilon/2)^{1/(k-1)}n^{(k-2)/(k-1)}\over 2n}<{1\over 4}
\]
Now observe that if $m_i\in[-\ell_n-1,\ell_n]$ and
\begin{eqnarray*}
{m_i(\epsilon/2)^{1/(k-1)}n^{(k-2)/(k-1)}\over 2n}&\le&{g_n((a_1-a_i)j)
\over n}\\
&\le &{(m_i+1)(\epsilon/2)^{1/(k-1)}n^{(k-2)/(k-1)}\over 2n},
\end{eqnarray*}
then
\[
\cos(2\pi(a_{i_1}-a_{i_2})j/n)\le\max\left(0.99,1-c_1{(\min(|m_i|,|m_i+1|))^2
(\epsilon/2)^{2/(k-1)}\over 8n^{2/(k-1)}}\right).
\]
Thus
\begin{eqnarray*}
\lefteqn{E(|\hat P_{a_1,\dots,a_k}(j)|^{2m})}\hspace{0.5in}\\
&\le&
b_1^m+\sum_{m_i\in[-\ell_n-1,\ell_n],i=2,\dots,k}{1.1(\epsilon/2)\over
2^{k-1}n}\\
&&\hspace{0.2in}\times\left(1-c_1(\min_{i_1\ne i_2}p_{i_1}p_{i_2})
\sum_{i=2}^k{(\min(|m_i|,|m_i+1|))^2(\epsilon/2)^{2/(k-1)}\over 
4n^{2/(k-1)}}\right)^m
\\
&=&b_1^m+2^{k-1}\sum_{m_i\in[0,\ell_n], i=2,\dots,k}{1.1(\epsilon/2)
\over 2^{k-1}n}\\
&&\hspace{0.2in}\times\left(1-c_1(\min_{i_1\ne i_2}p_{i_1}p_{i_2})
\sum_{i=2}^k{m_i^2(\epsilon/2)^{2/(k-1)}\over 4n^{2/(k-1)}}\right)^m.
\end{eqnarray*}

Since
\[
{m_i^2(\epsilon/2)^{2/(k-1)}\over 4n^{2/(k-1)}}<{1\over 16}\ ,
\]
$c_1\le 1$, and $\min_{i_1\ne i_2}p_{i_1}p_{i_2}\le 1/k$, we may 
conclude that
\[
c_1\left(\min_{i_1\ne i_2}p_{i_1}p_{i_2}\right)
\sum_{i=2}^k {m_i^2(\epsilon/2)^{2/(k-1)}\over 4n^{2/(k-1)}}<1.
\]
Thus for some constant $c_2>0$, 
\begin{eqnarray*}
\lefteqn{\left(1-c_1\left(\min_{i_1\ne i_2}p_{i_1}p_{i_2}\right)
\sum_{i=2}^k{m_i^2(\epsilon/2)^{2/(k-1)}\over 4n^{2/(k-1)}}\right)^m}
\hspace{0.5in}\\
&\le&\exp\left(-mc_2\sum_{i=2}^k{m_i^2(\epsilon/2)^{2/(k-1)}
\over 4n^{2/(k-1)}}\right).
\end{eqnarray*}
Thus
\begin{eqnarray*}
\lefteqn{E(|\hat P_{a_1,\dots,a_k}(j)|^{2m})}\hspace{0.5in}\\
&\le&b_1^m+\sum_{m_i\in[0,\ell_n], i=2,\dots,k}1.1{\epsilon \over 2n}
\exp\left(-mc_2\sum_{i=2}^km_i^2{(\epsilon/2)^{2/(k-1)}\over
4n^{2/(k-1)}}\right)\\
&\le&b_1^m+\sum_{m_i\in[0,\infty), i=2,\dots,k}1.1{\epsilon \over 2n}
\exp\left(-mc_2\sum_{i=2}^km_i{(\epsilon/2)^{2/(k-1)}\over 4n^{2/(k-1)}}
\right)\\
&=&b_1^m+\left({1.1\epsilon\over 2n}\right)/\left(1-\exp\left(-mc_2
{(\epsilon/2)^{2/(k-1)}\over 4n^{2/(k-1)}}\right)\right)^{k-1}.
\end{eqnarray*}
For some constant $\gamma>0$ and $m=\lfloor \gamma n^{2/(k-1)}\rfloor$,
\[
\lim_{n\rightarrow\infty}\left(1-\exp\left(-mc_2{(\epsilon/2)^{2/(k-1)}
\over 4n^{2/(k-1)}}\right)\right)^{k-1}\ge 0.7.
\]
Thus
\[
E(|\hat P_{a_1,\dots,a_k}(j)|^{2m})\le b_1^m+{0.9\epsilon\over n}
<{\epsilon\over n}
\]
for sufficiently large $n$, and the theorem follows from the Upper Bound Lemma.
\qed

\subsection{Random random walks on $\grpz_n$ where $n$ is not prime}

Dai and Hildebrand~\cite{daihild} generalized the result of 
Theorem~\ref{mvhptrfthm} to the case where $n$ need not be prime.
One needs to be cautious in determining which values to pick
for $a_1,\dots,a_k$. For example, if $n$ is even and $a_1,\dots,a_k$
are all even, then $X_m$ can never be odd. For another example, if $n$
is even and $a_1,\dots,a_k$ are all odd, then $X_m$ is never odd if $m$ is 
even and $X_m$ is never even if $m$ is odd. In both cases, 
$\|P_{a_1,\dots,a_k}^{*m}-U\|\not\rightarrow 0$ as $m\rightarrow\infty$.
Furthermore, if you choose $(a_1,\dots,a_k)$ uniformly from $G^k$, there
is a $1/2^k$ probability that $a_1,\dots,a_k$ are all even and
a $1/2^k$ probability that $a_1,\dots,a_k$ are all odd.

The following exercises develop a useful condition.

{EXERCISE.} Let $a_1,\dots,a_k\in\{0,\dots,n-1\}$.
Prove that the subgroup of $\grpz_n$ generated by
$\{a_2-a_1, a_3-a_1, \dots, a_k-a_1\}$ is $\grpz_n$ if and only if
$(a_2-a_1, a_3-a_1, \dots, a_k-a_1, n)=1$ where $(b_1, \dots, b_{\ell})$
is the greatest common divisor of $b_1, \dots, b_{\ell}$ in $\grpz$.

{EXERCISE.} Prove that $\|P_{a_1,\dots,a_k}^{*m}-U\|\rightarrow
0$ as $m\rightarrow\infty$ if and only if the subgroup of $\grpz_n$
generated by $\{a_2-a_1,\dots,a_k-a_1\}$ is $\grpz_n$. (Hint:
$X_m$ equals $m$ times $a_1$ plus some element of this subgroup.)

The main result of \cite{daihild} is
\begin{theorem}
\label{daihildthm}
Let $k\ge 2$ be a constant integer. Choose the set $S:=\{a_1,\dots,a_k\}$
uniformly from all  subsets of size $k$ from $\grpz_n$ such that
$(a_2-a_1, \dots, a_k-a_1, n)=1$ and $a_1, \dots, a_k$ are all
distinct. Suppose $p_1, \dots, p_k$ are positive constants with
$\sum_{i=1}^kp_i=1$. Then $E(\|P_{a_1,\dots,a_k}^{*m}-U\|)\rightarrow 0$
as $n\rightarrow\infty$ where $m:=m(n)\ge\sigma(n)n^{2/(k-1)}$ and
$\sigma(n)$ is any function with $\sigma(n)\rightarrow\infty$ as
$n\rightarrow\infty$. The expected value comes from the choice of the
set $S$.
\end{theorem}

The case $k=2$ can be handled in a manner similar to the case $k=2$ for
$n$ being prime. So we assume $k\ge 3$. Also, we shall let
$\epsilon>0$ be such that $1/(k-1)+(k+1)\epsilon<1$ throughout the
proof of Theorem~\ref{daihildthm}. As in the proof of Theorem~\ref{mvhptrfthm},
we can use the Upper Bound Lemma to bound $E(\|P_{a_1,\dots,a_k}^{*m}-U\|^2)$.
We may write $P$ instead of $P_{a_1,\dots,a_k}$.

We shall consider $3$ categories of values for $j$.
The proofs for the first and third categories use ideas from \cite{daihild}.
The proof for the second category at times diverges from the proof in
\cite{daihild}; the ideas in this alternate presentation were discussed in
personal communications between the authors of \cite{daihild}.

The first category has $J_1:=\{j:(j,n)>n^{(k-2)/(k-1)+\epsilon}$ and 
$1\le j\le n-1\}$. In this case, we have the following lemma.
\begin{lemma}
\label{largegcd}
$\sum_{j\in J_1}|\hat P(j)|^{2m}\rightarrow 0$ as $n\rightarrow\infty$ for 
$m\ge \sigma(n)n^{2/(k-1)}$ where $\sigma(n)\rightarrow\infty$ as $n\rightarrow
\infty$.
For a given $\sigma(n)$, this convergence is uniform over all choices of
the set $S$ where $a_1,\dots,a_k$ are distinct and
$(a_2-a_1,\dots,a_k-a_1,n)=1$.
\end{lemma}
{\bf Proof of Lemma:} Observe that if $n$ divides $ja_h-ja_1$ for all
$h=2,\dots,k$, then $n$ would divide $j$ since $(a_2-a_1,\dots,a_k-a_1,n)=1$.
So for some value $h\in\{2,\dots,k\}$, $n$ does not divide
$ja_h-ja_1$. Let $\omega=e^{2\pi i/n}$. Thus, with this
value of $h$, we get
\begin{eqnarray*}
|\hat P(j)|&=&|p_1\omega^{ja_1}+\dots+p_k\omega^{ja_k}|\\
&\le&1-2\min(p_1,p_h)+\min(p_1,p_h)|\omega^{ja_1}+\omega^{ja_h}|\\
&=&1-2\min(p_1,p_h)(1-|\cos(\pi j(a_h-a_1)/n)|)\\
&\le&1-ca^2/n^2
\end{eqnarray*}
where $a:=(j,n)$ and $c>0$ is a constant not depending on $S$. To see the last
inequality, observe that since $j(a_h-a_1)\not\equiv 0 \pmod n$, we get
$(j/a)(a_h-a_1)\not\equiv 0\pmod{n/a}$ and $(j/a)(a_h-a_1)\in\grpz$.
Thus $|\cos(\pi j(a_h-a_1)/n)|\le\cos(\pi a/n)$.

The proof of the lemma is completed with the following exercise.

{EXERCISE.} Show that for $m$ in the lemma and with 
$a>n^{(k-2)/(k-1)+\epsilon}$, we get $(1-ca^2/n^2)^m<c_1\exp(-n^{2\epsilon})$
for some constant $c_1>0$ and $n(1-ca^2/n^2)^{2m}\rightarrow 0$
as $n\rightarrow\infty$.

\qed

The next category of values of $j$ is $J_2:=\{j:(j,n)<n^{(k-2)/(k-1)-\epsilon}$
and $1\le j\le n-1\}$.
\begin{lemma}
\label{smallgcd}
$\sum_{j\in J_2}E(|\hat P(j)|^{2m})\rightarrow 0$ if $m\ge \sigma(n)n^{2/(k-1)}$
where $\sigma(n)\rightarrow\infty$ as $n\rightarrow\infty$ where the 
expectation is as in Theorem~\ref{daihildthm}.
\end{lemma}

To prove this lemma, we use the following.
\begin{lemma}
\label{smallgcdalt}
$\sum_{j\in J_2}E(|\hat P(j)|^{2m})\rightarrow 0$ if $m\ge \sigma(n)n^{2/(k-1)}$
where $\sigma(n)\rightarrow\infty$ as $n\rightarrow\infty$ where the 
expectation is over a uniform choice of $(a_1,\dots,a_k)\in G^k$.
\end{lemma}

{\bf Proof of Lemma~\ref{smallgcdalt}:} Although Dai and 
Hildebrand~\cite{daihild} used a different method, the proof of this lemma
can proceed in a manner similar to the proof of Theorem~\ref{mvhptrfthm}
for $k\ge 3$; however, one must be careful in proving the analogue of
Lemma~\ref{cubelemma}. Note in particular that instead of ranging over
$-(n-1)/2,\dots,-2,-1,1,2,\dots,(n-1)/2$, the value $g_n((a_1-a_i)j)$
ranges over multiples of $(j,n)$ in $(-n/2,n/2]$. Furthermore, distinct values
of $a_i$ need not lead to distinct values of $g_n((a_1-a_i)j)$.

{EXERCISE.} Show that Lemma~\ref{cubelemma} still holds if $n$ is not
prime but $j\in J_2$ and $(a_1,\dots,a_k)$ are chosen uniformly from
$G^k$. Then show Lemma~\ref{smallgcdalt} holds.

\qed

To complete the proof of Lemma~\ref{smallgcd}, note
that if $d$ is a divisor of $n$ and $(a_1,\dots,a_k)$ is
chosen uniformly from $G^k$, then the probability that $d$ divides all
of $a_2-a_1, \dots, a_k-a_1$ is $1/d^{k-1}$. Thus, given $n$, the probability
that $a_2-a_1, \dots, a_k-a_1$ have a common divisor greater than $1$ is
at most $\sum_{d=2}^{\infty}d^{-(k-1)}\le(\pi^2/6)-1<1$ if $k\ge 3$.
Also observe that the probability that a duplication exists on $a_1,\dots,a_k$
approaches to $0$ as $n\rightarrow\infty$.

{EXERCISE.} Show that these conditions and Lemma~\ref{smallgcdalt} imply
Lemma~\ref{smallgcd}.

\qed

The last category of values for $j$ is $J_3:=\{j: n^{(k-2)/(k-1)-\epsilon}
\le (j,n)\le n^{(k-2)/(k-1)+\epsilon}$ and $1\le j\le n-1\}$.
\begin{lemma}
\label{intermediategcd}
$\sum_{j\in J_3}E(|\hat P(j)|^{2m})\rightarrow 0$ if $m\ge \sigma(n)
n^{2/(k-1)}$ with $\sigma(n)\rightarrow \infty$ as $n\rightarrow \infty$
and with the expectation as in Theorem~\ref{daihildthm}.
\end{lemma}

{\bf Proof:}
Since $\|P_{a_1,a_2,\dots,a_k}^{*m}-U\|=\|P_{0,a_2-a_1,\dots,a_k-a_1}^{*m}-U\|$,
we may assume without loss of generality that $a_1=0$. Hence we assume
$(a_2,\dots,a_k,n)=1$. Let
$\langle x\rangle$ denote the distance of $x$ from the nearest multiple of $n$.
If $\langle ja_{\ell}\rangle>n^{(k-2)/(k-1)+\epsilon}$, then by reasoning
similar to that in the proof of Lemma~\ref{largegcd}, we may conclude that
$|\hat P(j)|^{2m}<c_1\exp(-n^{2\epsilon})$ for some constant $c_1>0$.
Otherwise $\langle ja_{\ell}\rangle\le n^{(k-2)/(k-1)+\epsilon}$ for
$\ell=2,\dots k$; let $B$ be the number of $(k-1)$-tuples $(a_2,\dots,a_k)$
satisfying this condition. Since $(j,n)\le n^{(k-2)/(k-1)+\epsilon}$, we may
conclude by Proposition~\ref{numbervalues} below that for some positive 
constant $c_2$, we have $B<c_2(n^{(k-2)/(k-1)+\epsilon})^{k-1}$.
Thus 
\[
E(|\hat P(j)|^{2m})<c_1\exp(-n^{2\epsilon})+{B\over Tn^{k-1}b(n)}
\]
where $T:=1-\sum_{d=2}^{\infty}d^{-(k-1)}$ and $b(n)$ is the probability
that $(a_2,\dots,a_k)$ when chosen at random from $G^{k-1}$
has no $0$ coordinates and no pair of coordinates 
with the same value.
Note that $b(n)\rightarrow 1$ as $n\rightarrow\infty$. Thus $E(|\hat P(j)|^{2m})
<c_1\exp(-n^{2\epsilon})+c_3n^{(k-1)\epsilon-1}$ for some constant $c_3>0$.

Let $D$ be the number of divisors of $n$. By Proposition~\ref{numberdivisors}
below, $D<c_4n^{\epsilon}$ for some positive constant $c_4$. Also, if $a$
divides $n$, then there are at most $n/a$ natural numbers in $[1,n-1]$
with $(j,n)=a$. For $j\in J_3$,
\[
{n\over a}\le{n \over n^{(k-2)/(k-1)-\epsilon}}=n^{1/(k-1)+\epsilon}.
\]
Thus
\[
\sum_{j\in J_3}E(|\hat P(j)|^{2m})\le c_4n^{\epsilon}n^{1/(k-1)+\epsilon}
(c_1\exp(-n^{2\epsilon})+c_3n^{(k-1)\epsilon-1})\rightarrow 0
\]
as $n\rightarrow \infty$ if $m>n^{2/(k-1)}\sigma(n)$. Note that we
used $1/(k-1)+(k+1)\epsilon<1$ here.   \qed

We need two propositions mentioned above.

\begin{proposition}
\label{numbervalues}
If $(j,n)\le b$ with $j\in\{1,\dots,n-1\}$, 
then the number of values $a$ in $0,1,\dots,n-1$ such that
$\langle ja\rangle\le b$ is less than or equal to $3b$.
\end{proposition}
The proof of this proposition may be found in \cite{daihild} or may be
done as an exercise.

\begin{proposition}
\label{numberdivisors}
For any $\epsilon$ with $0<\epsilon<1$, there is a positive constant 
$c=c(\epsilon)$ such that $d(n)\le cn^{\epsilon}$ for any natural number $n$
where $d(n)$ is the number of divisors of $n$.
\end{proposition}

{\bf Proof:} Suppose $n=p_1^{a_1}\dots p_r^{a_r}$ where $p_1, \dots, p_r$
are distinct prime numbers and $a_1,\dots,a_r$ are positive integers. Then
$d(n)=(a_1+1)(a_2+1)\dots(a_r+1)$. Note that $a_1, a_2, \dots$ here are
not the values selected from $\grpz_n$ and $p_1, p_2, \dots$ are not
probabilities! Let $M=e^{1/\epsilon}$. If $p_i>M$, then 
$(p_i^{a_i})^{\epsilon}\ge ((e^{1/\epsilon})^{a_i})^{\epsilon}=e^{a_i}
\ge 1+a_i$. If $p_i\le M$, then $(p_i^{a_i})^{\epsilon}\ge 2^{a_i\epsilon}=
e^{\epsilon a_i\log 2}\ge 1+\epsilon a_i\log 2\ge \epsilon(\log 2)(1+a_i)$.
Thus
\begin{eqnarray*}
n^{\epsilon}&=&
\left(\prod_{i=1}^rp_i^{a_i}\right)^{\epsilon}\\
&=&\left(\prod_{p_i\le M}p_i^{a_i}\right)^{\epsilon}
 \times \left(\prod_{p_i>m}p_i^{a_i}\right)^{\epsilon}\\
&\ge&
\prod_{p_i\le M}\epsilon(\log 2)(1+a_i) \times \prod_{p_i>M}(1+a_i)\\
&\ge&(\epsilon\log 2)^M\prod_{i=1}^r(1+a_i)
\end{eqnarray*}
since $\epsilon\log 2<1$. Thus $n^{\epsilon}\ge(\epsilon\log 2)^Md(n)$
and $d(n)\le cn^{\epsilon}$ where $c=(\epsilon\log 2)^{-M}$.\qed

{EXERCISE.} Give a rough idea of what $c$ is if $\epsilon=0.1$. 
(Hint: It's ridiculously large!)

{PROBLEM FOR FURTHER STUDY.} Give an argument so that 
Lemma~\ref{intermediategcd} or Theorem~\ref{daihildthm} does not need such a 
ridiculously large constant.

\subsection{Dou's version of the Upper Bound Lemma}

The Upper Bound Lemma of Diaconis and Shahshahani is particularly useful
for random walks on abelian groups or random walks where $P$ is constant
on conjugacy classes of $G$. Dou~\cite{dou} has adapted this lemma to a form
which is useful for some random random walks. This form was used in \cite{dou} 
to study some random random walks on various abelian groups, and Dou and
Hildebrand~\cite{douhild} extended some results in \cite{dou} involving
random random walks on abelian groups. Dou's lemma is the following.

\begin{lemma}
\label{doulemma}
Let $Q$ be a probability on a group $G$ of order $n$. Then for any
positive integer $m$,
\[
4\|Q^{*m}-U\|^2\le\sum_{\Omega}nQ(x_1)\dots Q(x_{2m})-\sum_{G^{2m}}
Q(x_1)\dots Q(x_{2m})
\]
where $G^{2m}$ is the set of all $2m$-tuples $(x_1,\dots,x_{2m})$ with 
$x_i\in G$ and $\Omega$ is the subset of $G^{2m}$ consisting of all $2m$-tuples
such that $x_1x_2\dots x_m=x_{m+1}x_{m+2}\dots x_{2m}$.
\end{lemma}

{\bf Proof:} The proof presented here uses arguments given 
in \cite{dou} and \cite{douhild}.
Label all the non-equivalent irreducible representations of $G$
by $\rho_1,\dots,\rho_h$. Assume $\rho_h$ is the trivial representation
and the representations are all unitary. Let $\chi_1,\dots,\chi_h$ be
the corresponding characters and $d_1,\dots,d_h$ be the corresponding
degrees. Note that $\rho_i(x)^{*}=(\rho_i(x))^{-1}=\rho_i(x^{-1})$ for
all $x\in G$ since $\rho_i$ is unitary. Thus
\[
\hat Q(\rho_i)^m=\sum_{x_1,\dots,x_m}Q(x_1)\dots Q(x_m)\rho_i(x_1\dots x_m)
\]
and
\[
(\hat Q(\rho_i)^m)^{*}=\sum_{x_{m+1},\dots,x_{2m}}Q(x_{m+1})\dots Q(x_{2m})
\rho_i((x_{m+1}\dots x_{2m})^{-1}))
\]
Thus if $s=(x_1\dots x_n)(x_{m+1}\dots x_{2m})^{-1}$, we get
\begin{eqnarray*}
\sum_{i=1}^{h-1}d_i Tr(\hat Q(\rho_i)^m(\hat Q(\rho_i)^m)^{*})
&=&\sum_{i=1}^{h-1}\sum_{G^{2m}}Q(x_1)\dots Q(x_{2m})d_i\chi_i(s)\\
&=&\sum_{G^{2m}}Q(x_1)\dots Q(x_{2m})\sum_{i=1}^{h-1}d_i\chi_i(s)
\end{eqnarray*}
Note that $d_h\chi_h(s)=1$ for all $s\in G$ and
\[
\sum_{i=1}^h d_i\chi_i(s)=\cases{n&if $s=e$\cr 0&otherwise}
\]
where $e$ is the identity element of $G$. The lemma follows from the
Upper Bound Lemma. \qed

To use this lemma, we follow some notation as in \cite{dou} and \cite{douhild}.
We say a
$2m$-tuple $(x_1,\dots,x_{2m})$
is of {\it size} $i$ if the set $\{x_1,\dots,x_{2m}\}$ has
$i$ distinct elements. An {\it $i$-partition} of $\{1,\dots,2m\}$ is a set
of $i$ disjoint subsets $\tau=\{\Delta_1,\dots,\Delta_i\}$ such that
$\Delta_1\cup\dots\cup\Delta_i=\{1,\dots,2m\}$. An {\it $i$-partition}
of the number $2m$ is an $i$-tuple of integers $\pi=(p_1,\dots,p_i)$ such
that $p_1\ge \dots\ge p_i\ge 1$ and $\sum_{j=1}^ip_j=2m$.

Note that each $2m$-tuple in $G^{2m}$ of size $i$ gives rise to an
$i$-partition of $2m$ in a natural way.
For $1\le j\le i$, let $\Delta_j\subset\{1,\dots,2m\}$ be a maximal
subset of indices for which the corresponding coordinates are the same.
Then $\Delta_1,\dots,\Delta_i$ form an $i$-partition of $\{1,\dots,2m\}$;
we call this $i$-partition the {\it type} of the $2m$-tuple. If
$|\Delta_1|\ge\dots\ge|\Delta_i|$, then $\pi=(|\Delta_1|,\dots,|\Delta_i|)$
is an $i$-partition of $2m$, and we say the type $\tau$ corresponds to
$\pi$.

{EXAMPLE.} If $\nu=(0,11,5,5,1,3,0,5)\in\grpz_{12}^8$, then the type
of $\nu$ is $\tau=\{\{3,4,8\},\{1,7\},\{2\},\{5\},\{6\}\}$ and the
corresponding $5$-partition of the number $8$ is $\pi=(3,2,1,1,1)$.

Now suppose that $\tau=\{\Delta_1,\dots,\Delta_i\}$
is a type corresponding to a partition $\pi$ of $2m$. Let $N_{\pi}(\tau)$ be
the number of $2m$-tuples in $\Omega$ of type $\tau$. 

A little thought should give the following lemma.
\begin{lemma} 
\label{doucount}
$N_{\pi}(\tau)$ is the number of $i$-tuples $(y_1,\dots,y_i)$ with
distinct coordinates in $G$ that are solutions to the induced equation obtained
from $x_1\dots x_m=x_{m+1}\dots x_{2m}$ by substituting $y_j$ for $x_{\ell}$
if $\ell\in\Delta_j$.
\end{lemma}

{EXAMPLE.} If $\tau=\{\{3,4,8\},\{1,7\},\{2\},\{5\},\{6\}\}$, then this 
equation is
$y_2y_3y_1^2=y_4y_5y_2y_1$.

{EXERCISE.} Suppose $\nu=(0,3,4,1,5,1,4,5,3,3)\in\grpz_{12}^{10}$.
Find the type $\tau$ of $\nu$. Then find the 
value $i$ and the induced equation described in the previous lemma.

We can adapt Lemma~\ref{doulemma} to prove the following lemma. The lemma and
its proof are minor modification of some material presented in
\cite{dou} and \cite{douhild}. In particular, the argument in \cite{dou}
covered a broader range of probabilities $Q$.

\begin{lemma}
\label{doupart}
Suppose $(a_1,\dots,a_k)$ are chosen uniformly from all $k$-tuples with
distinct elements of $G$. Also suppose that $Q(a_i)=1/k$. Then
\[
E(\|Q^{*m}-U\|^2)\le\sum_{i=1}^{\min(k,2m)}\sum_{\pi\in P(i)} 
{1\over k^{2m}}{[k]_i\over [n]_i}\sum_{\tau\in T(\pi)}(nN_{\pi}(\tau)-[n]_i)
\]
where $[n]_i=n(n-1)\dots(n-i+1)$, $P(i)$ is the set of all $i$-partitions
of $2m$, and $T(\pi)$ is the set of all types which correspond to $\pi$.
\end{lemma}

{\bf Proof:} Suppose $\pi$ is an $i$-partition of $2m$. A $2m$-tuple of
$\pi$ is defined to be a $2m$-tuple whose type corresponds to $\pi$, Let
$D_1(\pi)$ be the set of all $2m$-tuples of $\pi$ in $\Omega$, and let
$D_2(\pi)$ be the set of all $2m$-tuples of $\pi$ in $G^{2m}$. Thus
$|D_1(\pi)|=\sum_{\tau\in T(\pi)}N_{\pi}(\tau)$ and
$|D_2(\pi)|=\sum_{\tau\in T(\pi)}M_{\pi}(\tau)$ where
$M_{\pi}(\tau)$ is the number of $2m$-tuples of type $\tau$ in $G^{2m}$
where $\pi$ is the corresponding $i$-partition of the number $2m$. It can
readily be shown that $M_{\pi}(\tau)=[n]_i$.

Observe from Lemma~\ref{doulemma} that
\begin{eqnarray*}
4E(\|Q^{*m}-U\|^2)
&\le&\sum_{i=1}^{2m}\sum_{\pi\in P(i)}
\left(\sum_{(x_1,\dots,x_{2m})
\in D_1(\pi)}nE(Q(x_1)\dots Q(x_{2m}))\right.\\
&&\hspace{0.25in} \left.-\sum_{(x_1,\dots,x_{2m})\in D_2(\pi)}
E(Q(x_1)\dots Q(x_{2m}))\right).
\end{eqnarray*}

Now let's consider $E(Q(x_1)\dots Q(x_{2m}))$. This expectation depends only
on the size $i$ of $(x_1,\dots,x_{2m})$. The probability that a given 
$i$-tuple $(y_1,\dots,y_i)$ with distinct elements of $G$ is contained in
a random sample of $G$ is $[k]_i/[n]_i$. Thus
\[
E(Q(x_1)\dots Q(x_{2m}))={1\over k^{2m}}{[k]_i\over [n]_i}
\]
Note that if the size of $(x_1,\dots,x_{2m})$ is greater than $k$, then
$Q(x_1)\dots Q(x_{2m})$ must be $0$.

Thus
\[
\sum_{x\in D_1(\pi)}nE(Q(x_1)\dots Q(x_{2m}))={1\over k^{2m}}{[k]_i\over
[n]_i}\sum_{\tau\in T(\pi)}nN_{\pi}(\tau)
\]
and
\begin{eqnarray*}
\sum_{x\in D_2(\pi)}E(Q(x_1)\dots Q(x_{2m}))&=&{1\over k^{2m}}{[k]_i\over
[n]_i}\sum_{\tau\in T(\pi)}M_{\pi}(\tau)\\
&=&{1\over k^{2m}}{[k]_i\over [n]_i}\sum_{\tau\in T(\pi)}[n]_i
\end{eqnarray*}
The lemma follows by substitution and easy algebra.
\qed

Next we consider a result in Dou~\cite{dou}. The next theorem 
and its proof essentially come from \cite{dou} but use a simpler and less 
general expression for the probability $P$.

We assume that $G$ is an abelian group with $n$ elements such that 
$n=n_1\dots n_t$ where $n_1\ge \dots\ge n_t$ are prime numbers, $t\le L$ for
some value $L$ not depending on $n$, and $n_1\le An_t$ for some value
$A$ not depending on $n$.

\begin{theorem}
\label{douabelian}
Suppose $G$ satisfies the conditions in the previous paragraph and
$k>2L+1$ is constant. Suppose $(a_1,\dots,a_k)$ is chosen uniformly 
from $k$-tuples with distinct elements of $G$. Then
for some function $f(n)\rightarrow 0$ as $n\rightarrow\infty$
(with $f(n)$ not depending on the choice of $G$)
\[
E(\|P_{a_1,\dots,a_k}^{*m}-U\|)\le f(n)
\]
where $m=c(n)n^{2/(k-1)}$ where $c(n)\rightarrow
\infty$ as $n\rightarrow\infty$ and $p_1=\dots=p_k=1/k$ so that
$P_{a_1,\dots,a_k}(s)=1/k$ if $s=a_i$ for some $i$ in $1,\dots,k$.
\end{theorem}

{\bf Proof:} Without loss of generality, we may assume that 
\[
c(n)<n^{(1/L)-(2/(k-1))}A^{-1+(1/L)}.
\]
We may use Lemma~\ref{doupart}. Let
\[
B_1=\sum_{i=1}^{k-1}\sum_{\pi\in P(i)}{1\over k^{2m}}{[k]_i\over
[n]_i}\sum_{\tau\in T(\pi)}(nN_{\pi}(\tau)-[n]_i).
\]
It can be readily shown from Lemma~\ref{doucount} that
$nN_{\pi}(\tau)-[n]_i\le n[n]_i$. Thus
\begin{eqnarray*}
B_1&\le& \sum_{i=1}^{k-1}\sum_{\pi\in P(i)}{1\over k^{2m}}{[k]_i\over [n]_i}
\sum_{\tau\in T(\pi)}n[n]_i
\\
&=&{1\over k^{2m}}n\sum_{i=1}^{k-1}[k]_iS_{2m,i}
\end{eqnarray*}
where $S_{2m,i}$ is a Stirling number of the second kind; this
number is the number of ways to place $2m$ labeled balls in $i$
unlabeled boxes such that there are no empty boxes.

{EXERCISE.} Show that $(k-1)^{2m}=\sum_{i=1}^{k-1}[k-1]_iS_{2m,i}$. (Hint:
The left side is the number of ways to place $2m$ labeled balls in
$k-1$ labeled boxes where some boxes may be empty.)

Observe that $[k]_i=k[k-1]_i/(k-i)\le k[k-1]_i$ if $i\le k-1$. Thus
$B_1\le k(1/k^{2m})n(k-1)^{2m}\rightarrow 0$ as $n\rightarrow\infty$
for the specified $m$.

Now let 
\[
B_2=\sum_{\tau}{1\over k^{2m}}{[k]_k\over [n]_k}(nN_{\pi}(\tau)-[n]_k)
\]
where the sum is over all $k$-tuples $\tau$, i.e. the set of all
$k$-partitions of the set $\{1,2,\dots,2m\}$. We use the following
lemma (which, along with is proof, is based upon \cite{dou}).

\begin{lemma}
If $G$ is an abelian group satisfying the conditions for 
Theorem~\ref{douabelian} and $m=c(n)n^{2/(k-1)}$ where $c(n)\rightarrow
\infty$ as $n\rightarrow\infty$ such that 
\[
c(n)<n^{(1/L)-(2/(k-1))}A^{-1+(1/L)},
\]
then for each $k$-type $\tau$, either $N_{\pi}(\tau)=[n]_k$
or $N_{\pi}(\tau)\le [n]_{k-1}$. If $T_1=\{k-{\rm types\ }
\tau|N_{\pi}(\tau)=[n]_k\}$
and $T_2=\{k-{\rm types\ }\tau|N_{\pi}(\tau)\le[n]_{k-1}\}$, then
$|T_1|+|T_2|=S_{2m,k}$ and 
\[
|T_1|\le\kappa_{m,k}:=\sum_{r_1+\dots+r_k=m, r_1,\dots,r_k\ge 0}
{m\choose r_1,\dots,r_k}^2
\]
\end{lemma}

{\bf Proof:} We start with an exercise.

{EXERCISE.} Show that the conditions on $G$ and the restrictions on
$c(n)$ imply that $m<n_t$.

By Lemma~\ref{doucount}, $N_{\pi}(\tau)$ is the number of $k$-tuples
$(y_1,\dots,y_k)$ with distinct coordinates such that 
$\lambda_1y_1+\dots+\lambda_ky_k=0$ for some integers
$\lambda_1,\dots,\lambda_k$ with $|\lambda_i|\le m<n_t$.
Note that here we are using the fact that $G$ is abelian. Also note
 that $\lambda_i$
is the number of times $y_i$ is substituted for $x_j$ with $1\le j \le m$
minus the number of times $y_i$ is substituted for $x_j$ with
$m+1\le j\le 2m$; in other words, $\lambda_i=|\Delta_i\cap\{1,\dots,m\}|-
|\Delta_i\cap\{m+1,\dots,2m\}|$.
If $\lambda_1=\dots=\lambda_k=0$, then $N_{\pi}(\tau)=[n]_k$;
otherwise if $\lambda_j\ne 0$ for some $j$, then $y_j$ is solvable
in terms of the other variables in $y_1,\dots,y_k$ since $m<p_t$ and thus
$N_{\pi}(\tau)\le[n]_{k-1}$. Thus we may
conclude that $\lambda_1=\dots=\lambda_k=0$ if and only if
$\tau=\{\Delta_1^{\prime}\cup\Delta_1^{\prime\prime},\dots,
\Delta_k^{\prime}\cup\Delta^{\prime\prime}\}$ where
$\tau^{\prime}=\{\Delta_1^{\prime},\dots,\Delta_k^{\prime}\}$
and
$\tau^{\prime\prime}=\{\Delta_1^{\prime\prime},\dots,\Delta_k^{\prime\prime}\}$
are $k$-types of $\{1,\dots,m\}$ and $\{m+1,\dots,2m\}$ respectively with
$|\Delta_i^{\prime}|=|\Delta_i^{\prime\prime}|$ for $i=1,2,\dots,k$.
The inequality $|T_1|\le \kappa_{m,k}$ is elementary. The equality
$|T_1|+|T_2|=S_{2m,k}$ follows quickly from the definitions of $S_{2m,k}$,
$T_1$, and $T_2$. The lemma is thus proved. \qed

Thus $B_2\le B_{2,1}+B_{2,2}$ where
\[
B_{2,1}=\sum_{\tau\in T_1}{1\over k^{2m}}{[k]_k\over [n]_k}(n-1)[n]_k
\]
and
\[
B_{2,2}=\sum_{\tau\in T_2}{1\over k^{2m}}{[k]_k\over [n]_k}[n]_{k-1}(k-1).
\]
Note that in defining $B_{2,2}$, we used the fact that
$n[n]_{k-1}-[n]_k=[n]_{k-1}(k-1)$. Now observe that
\begin{eqnarray*}
B_{2,2}&\le &(k-1){1\over k^{2m}}[k]_k|T_2|/(n-k+1)\\
&\le &(k-1){1\over k^{2m}}[k]_kS_{2m,k}/(n-k+1)\\
&\le &(k-1){1\over k^{2m}}k^{2m}/(n-k+1)
\end{eqnarray*}
since $S_{2m,k}\le k^{2m}/k!$ because $k^{2m}$ is the no more than the
number of ways to place $2m$ labeled balls in $k$ labeled boxes where
no box is left empty. Thus $B_{2,2}\rightarrow 0$ as $n\rightarrow \infty$.

We also have
\begin{eqnarray*}
B_{2,1}&\le &{1\over k^{2m}}[k]_k(n-1)|T_1|\\
&\le &{1\over k^{2m}}[k]_k(n-1)\kappa_{m,k}\\
&\le &{1\over k^{2m}}[k]_k(n-1)c_0k^{2m}/m^{(k-1)/2}
\end{eqnarray*}
since an appendix of \cite{dou} proves $\kappa_{m,k}\le c_0k^{2m}/m^{(k-1)/2}$
for some positive constant $c_0$. Thus for some value $c$ which may depend
on $k$ but not $n$, we get $B_{2,1}\le c(n-1)/m^{(k-1)/2}$. The theorem
follows. \qed

{PROBLEM FOR FURTHER STUDY.} The restrictions on $G$ in 
Theorem~\ref{douabelian} seem excessive, but we used the fact that $m<n_t$
to bound $|T_1|$. Perhaps arguments can be made for a broader range of
$G$ to deal with some cases where we do not have this fact. Indeed, similar
squares of multinomial coefficients, along with some additional terms, appear
in the argument which Dai and Hildebrand~\cite{daihild} use to 
prove Lemma~\ref{smallgcdalt} described earlier; 
this proof might serve as a starting point for a possible proof of an
extension of Theorem~\ref{douabelian}.

It should be noted that Greenhalgh~\cite{greenhalghpreprint}
has also used arguments involving squares of such multinomial coefficients
to prove results similar to Theorem~\ref{mvhptrfthm}.

Further results using these techniques appear in Dou and 
Hildebrand~\cite{douhild}. In particular, the following two theorems
are shown there; their proofs will not be presented here.

\begin{theorem}
\label{douhilda}
Suppose $G$ is an arbitrary finite group of order $n$ and
$k=\lfloor(\log n)^a\rfloor$ where $a>1$ is constant. Let $\epsilon>0$ 
be given. Suppose $(a_1,\dots,a_k)$ is chosen uniformly from
$k$-tuples with distinct elements of $G$. Then for some
function $f(n)\rightarrow 0$ as $n\rightarrow\infty$ (with
$f(n)$ not depending on the choice of $G$),
\[
E(\|P_{a_1,\dots,a_k}^{*m}-U\|)\le f(n)
\]
if
\[
m=m(n)>{a\over a-1}{\log n\over\log k}(1+\epsilon).
\]
\end{theorem}

\begin{theorem}
\label{douhildb}
Suppose $G$ is an arbitrary finite group of order $n$. Also suppose
$k=\lfloor a\log n\rfloor$ and $m=\lfloor b\log n\rfloor$ where
$a$ and $b$ are constants with
$a>e^2$, $b<a/4$, and $b\log(eb/a)<-1$. Suppose $(a_1,\dots,a_k)$ is
chosen uniformly from $k$-tuples with distinct elements of $G$. Then for
some function $f(n)\rightarrow 0$ as $n\rightarrow\infty$ (with $f(n)$
not depending on the choice of $G$),
\[
E(\|P_{a_1,\dots,a_k}^{*m}-U\|)\le f(n).
\]
\end{theorem}

Roichman~\cite{roichman} uses spectral arguments to get results similar
to these theorems and to extend them to symmetric random walks where at
each step one multiplies either by $a_i$ or by $a_i^{-1}$ (with probability
$1/2$ each) where $i$ is chosen uniformly from $\{1,\dots,k\}$.

\subsection{Extensions of a result of Erd\"os and R\'enyi}

Some early results involving probability on finite groups appear in
an article of Erd\"os and R\'enyi~\cite{erdosrenyi}.
In it, they give the following theorem.
\begin{theorem}
\label{erdosrenyitheorem}
Suppose $k\ge 2\log_2 n+2\log_2 (1/\epsilon)+\log_2(1/\delta)$ and
$J=(a_1,\dots,a_k)$ is a random $k$-tuple of elements from an abelian
group $G$ of order $n$. If $b\in G$, let $V_k(b)$ be the number of
$(\epsilon_1,\dots,\epsilon_k)\in\{0,1\}^k$ such that $b=\epsilon_1a_1+
\dots+\epsilon_ka_k$. Then
\[
Pr\left(\max_{b\in G}\left|V_k(b)-{2^k\over n}\right|\le\epsilon{2^k\over n}
\right)>1-\delta.
\]
\end{theorem}

Near the end of their paper, Erd\"os and R\'enyi note that this theorem
can be generalized to non-abelian groups $G_n$ of order $n$ by counting
the number of ways an element $b$ can be written in the form
$b=a_{i_1}a_{i_2}\dots a_{i_r}$ where $1\le  i_1<i_2<\dots<i_r\le k$ and
$0\le r\le k$. They find it unnatural to assume that if $i<j$ that
$a_i$ would have to appear before $a_j$. They also assert
that if we ``do not consider only such products in which $i_1<i_2<\dots
<i_r$, then the situation changes completely. In this case the
structure of the group $G_n$ becomes relevant.''

Further results by Pak~\cite{pak} and Hildebrand~\cite{mvhjtp} built
upon this result in \cite{erdosrenyi} and its extension to non-abelian
groups to get results for
random ``lazy'' random walks on arbitrary finite groups. These
``lazy'' random walks are such that at each step, there's a large
probability that the walk stays at the same group element. 
These results, despite
the assertion of \cite{erdosrenyi}, do not depend on the structure of the
group $G$ and involve products of the form $a_{i_1}a_{i_2}\dots
a_{i_r}$ where $i_1, i_2, \dots, i_r$ need not be in increasing order and 
may be repeated.

For the next two theorems, which come from \cite{mvhjtp}, we use the following 
notation. If $J=(a_1,\dots,a_k)\in G^k$, then
\[
P(s)=P_J(s)={1\over 2k}|\{i: a_i=s\}|+{1\over 2}\delta_{\{s=e\}}
\]
where $\delta_{\{s=e\}}$ is $1$ if $s$ is the identity $e$ of $G$ and
$0$ otherwise.
Note that this expression is different from the expression
$P_{a_1,\dots,a_k}(s)$ in the rest of this article.

\begin{theorem}
\label{jtpone}
Suppose $a>1$ and $\epsilon>0$ are given. Let $k=\lceil a\log_2 n\rceil$.
Suppose $m=m(n)>(1+\epsilon)a\log (a/(a-1))\log_2 n$. Then for some
function $f_1(n)\rightarrow 0$ as $n\rightarrow\infty$ (where $f_1(n)$
does not depend on which group $G$ of order $n$ is being considered),
$E(\|P^{*m}-U\|)\le f_1(n)$ as $n\rightarrow\infty$ where $J$ is chosen
uniformly from $G^k$.
\end{theorem}

\begin{theorem}
\label{jtptwo}
Suppose $k=\log_2 n+f(n)$ where $f(n)\rightarrow \infty$ as $n\rightarrow\infty$
and $f(n)/\log_2 n\rightarrow 0$ as $n\rightarrow\infty$. Let $\epsilon>0$
be given. If $m=m(n)>(1+\epsilon)(\log_2 n)(\log(\log_2 n))$, then for some
function $f_2(n)\rightarrow 0$ as $n\rightarrow\infty$ (where $f_2(n)$
does not depend on which group $G$ of order $n$ is being considered),
$E(\|P^{*m}-U\|)\le f_2(n)$ where $J$ is chosen uniformly from $G^k$.
\end{theorem}

The proofs of Theorems~\ref{jtpone} and \ref{jtptwo} use the following 
variation of Theorem~\ref{erdosrenyitheorem}.
In it, note that $g^0$ is the identity element of $G$.
\begin{lemma}
\label{variationerdosrenyi}
Let $J=(a_1,\dots,a_k)$. Suppose $j\le k$. Let $Q_J(s)$ be the
probability that $s=a_1^{\epsilon_1}a_2^{\epsilon_2}\dots a_j^{\epsilon_j}$
where $\epsilon_1,\epsilon_2,\dots,\epsilon_j$ are i.i.d. uniform on $\{0,1\}$.
If $J$ is chosen uniformly from all $k$-tuples $(a_1,\dots,a_k)$ of
elements of $G$, then $Pr(\|Q_J-U\|\le\epsilon)\ge 1-\delta$ for each\
$j\ge \log_2n+2\log_2(1/\epsilon)+\log_2(1/\delta)$.
\end{lemma}
{\bf Proof:} This proof follows \cite{mvhjtp} and extends a proof in
\cite{erdosrenyi}.

Let $V_j(s)=2^jQ_J(s)$. Observe that
\begin{eqnarray*}
4\|Q_J-U\|^2&=&\left(\sum_{s\in G}\left|{V_j(s)\over 2^j}-{1\over n}
\right|\right)^2\\
&\le &n\sum_{s\in G}\left({V_j(s)\over 2^j}-{1\over n}\right)^2
\end{eqnarray*}
by the Cauchy-Schwarz inequality; this argument is very similar to
part of the proof of the Upper Bound Lemma described in Chapter 3 of
Diaconis~\cite{diaconis}.

Thus
\begin{eqnarray*}
Pr(2\|Q_J-U\|>\epsilon)&\le &Pr\left(n\sum_{s\in G}\left(
{V_j(s)\over 2^j}-{1\over n}\right)^2>\epsilon^2\right)\\
&=&Pr\left(\sum_{s\in G}\left(V_j(s)-{2^j\over n}\right)^2>
{\epsilon^2 2^{2j}\over n}\right).
\end{eqnarray*}

It can be shown (as on p. 130 of \cite{erdosrenyi} extended to non-abelian
groups) that
\[
E\left(\sum_{s\in G}\left(V_j(s)-{2^j\over n}\right)^2\right)=
2^j\left(1-{1\over n}\right)
\]
where the expectation comes from choosing $J$ uniformly from all
$k$-tuples $(a_1,\dots,a_k)$ of elements of $G$.

Thus by Markov's inequality, we get
\begin{eqnarray*}
Pr(\|Q_J-U\|>\epsilon)&\le &Pr(2\|Q_J-U\|>\epsilon)\\
&\le &{2^j(1-(1/n))\over \epsilon^2 2^{2j}/n}\\
&\le &{n \over \epsilon^2 2^j}\le\delta.
\end{eqnarray*}
\qed

We say that a family of probability distributions $R_J$ depending on
$J\in G^k$ is $(\alpha,\beta)$-{\it good} in variation distance if
$Pr(\|R_J-U\|>\alpha)\le\beta$ where the probability is over a 
uniform choice of all $k$-tuples for $J$. Thus Lemma~\ref{variationerdosrenyi}
shows that $a_1^{\epsilon_1}\dots a_j^{\epsilon_j}$ is 
$(\epsilon,\delta)$-good in variation distance if 
$j\ge\log_2 n+2\log_2(1/\epsilon)+\log_2(1/\delta)$.

Theorems \ref{jtpone} and \ref{jtptwo} look at the variation distance
from the uniform distribution of a probability distribution of 
\[
a_{i_1}^{\epsilon_1}a_{i_2}^{\epsilon_2}\dots a_{i_m}^{\epsilon_m}
\]
where $i_1,\dots,i_m$ are i.i.d. uniform on $\{1,\dots,k\}$,
$\epsilon_1,\dots,\epsilon_m$ are i.i.d. uniform on $\{0,1\}$,
and $(i_1,\dots,i_m)$ and $(\epsilon_1,\dots,\epsilon_m)$ are
independent. Using Lemma~\ref{variationerdosrenyi} to examine this
distribution requires considerable care.

First let's consider the case where $i_1,\dots,i_m$ are all given
and consist of at least $j$ distinct values. Suppose also that the
value $\epsilon_{\ell}$ is given if $i_{\ell}=i_{\ell^{\prime}}$
for some $\ell^{\prime}<\ell$ or if $\{i_1,\dots,i_{\ell-1}\}$ has
at least $j$ distinct values; in other words, $\epsilon_{\ell}$ is given
if $i_{\ell}$ appeared earlier in the $m$-tuple
$(i_1,\dots,i_m)$ or if $i_{\ell}$ is not among the first
$j$ distinct values in this $m$-tuple. We assume that the
remaining $j$ values from $\epsilon_1,\dots,\epsilon_m$ are i.i.d.
uniform on $\{0,1\}$. For example, if $j=5$, $k=7$, and $m=9$, such
an expression may look like
\[
a_3^{\epsilon_1}a_4^{\epsilon_2}a_2^{\epsilon_3}a_4^1 a_1^{\epsilon_5}
a_3^1 a_6^{\epsilon_7}a_7^0 a_2^1
\] 
where $\epsilon_1$, $\epsilon_2$, $\epsilon_3$, $\epsilon_5$, and
$\epsilon_7$ are i.i.d. uniform on $\{0,1\}$.

To use Lemma~\ref{variationerdosrenyi} to examine such expressions, we need
to consider the ``pulling through'' technique described in Pak~\cite{pak}
and subsequently in Hildebrand~\cite{mvhjtp}.
\begin{proposition}
\label{pullthrough}
Suppose $h=a_1^{\epsilon_1}\dots a_{\ell}^{\epsilon_{\ell}}x
a_{\ell+1}^{\epsilon_{\ell+1}}\dots a_j^{\epsilon_j}$ where
$\epsilon_1,\dots,\epsilon_j$ are each in $\{0,1\}$ and $x$ is a fixed
function of $a_1,\dots,a_{\ell}$. Then
\[
h=a_1^{\epsilon_1}\dots a_{\ell}^{\epsilon_{\ell}}\left(a_{\ell+1}^x\right)^{
\epsilon_{\ell+1}}\dots\left(a_j^x\right)^{\epsilon_j}x
\]
where $g^x:=xgx^{-1}$. Furthermore if $a_1,\dots,a_j$ are i.i.d. uniform on
$G$, then $a_1,\dots,a_{\ell},a_{\ell+1}^x,\dots,a_j^x$ are i.i.d. uniform 
on $G$.
\end{proposition}

{\bf Proof:} The alternate expression for $h$ can be readily verified. Note 
that since $x$ does not depend on $a_{\ell+1}$ and since
$a_1,\dots,a_{\ell},a_{\ell+1}$ are i.i.d. uniform on $G$, the expression
$xa_{\ell+1}x^{-1}$ will be uniform on $G$ independent of $a_1,\dots,a_{\ell}$.
Continuing in the same way completes the proof of the proposition. \qed

This proposition can be used repeatedly. For example, if
\[
h=a_3^{\epsilon_1}a_4^{\epsilon_2}a_2^{\epsilon_3}a_4^1 a_1^{\epsilon_5}
a_3^1 a_6^{\epsilon_7}a_7^0 a_2^1,
\]
then
\[
h=a_3^{\epsilon_1}a_4^{\epsilon_2}a_2^{\epsilon_3}\left(a_1^{a_4}\right)^{
\epsilon_5}\left(a_6^{a_4a_3}\right)^{\epsilon_7}a_4a_3a_7^0 a_2^1.
\]
Furthermore, if $a_3$, $a_4$, $a_2$, $a_1$, and $a_6$ are i.i.d. uniform
on $G$, then so are $a_3$, $a_4$, $a_2$, $a_1^{a_4}$, and $a_6^{a_4a_3}$.
Also note that if $a_1, \dots, a_7$ are given and $\epsilon_1$, $\epsilon_2$,
$\epsilon_3$, $\epsilon_5$, and $\epsilon_7$ are i.i.d. uniform on $\{0,1\}$,
then the probabilities $P$ and $Q$ given by $P(s)=Pr(s=a_3^{\epsilon_1}
a_4^{\epsilon_2}a_2^{\epsilon_3}(a_1^{a_4})^{\epsilon_5}
(a_6^{a_4a_3})^{\epsilon_7})$ and $Q(s)=Pr(s=h)$ have the same
variation distance from the uniform distribution.

Thus we may conclude the following.
\begin{lemma}
\label{goodlemma}
Suppose $I=(i_1,\dots,i_m)$ where $i_1,\dots,i_m\in\{1,\dots,k\}$.
Suppose $I$ has at least $j$ distinct values where $j\ge \log_2n+
2\log_2(1/\epsilon)+\log_2(1/\delta)$ and $j\le k$. Suppose $\vec\epsilon$
is a vector determining $\epsilon_{\ell}$ if $i_{\ell}=i_{\ell^{\prime}}$
for some $\ell^{\prime}<\ell$ or if $\{i_1,\dots,i_{\ell-1}\}$ has at least
$j$ distinct values. Suppose the remaining $j$ values from $\epsilon_1,\dots,
\epsilon_m$ are i.i.d. uniform on $\{0,1\}$. Then $a_{i_1}^{\epsilon_1}
\dots a_{i_m}^{\epsilon_m}$ is $(\epsilon,\delta)$-good in variation
distance.
\end{lemma}

We need to put together probabilities for the various possibilities for $I$
and $\vec\epsilon$. The following exercise will be useful.

{EXERCISE. Suppose $P=p_1P_1+\dots+p_{\ell}P_{\ell}$ where $p_1,\dots,p_{\ell}$
are positive numbers which sum to $1$ and $P_1,\dots,P_{\ell}$ are 
probabilities on $G$. Show that
\[
\|P-U\|\le\sum_{j=1}^{\ell}p_j\|P_j-U\|.
\]
}

The following lemma comes from \cite{mvhjtp}.
\begin{lemma}
\label{puttogether}
Let $c>1$ be given. Suppose $j$ is given such that $j\le k$,
$j\ge\log_2 n+2\log_2(1/\alpha)+\log_2(1/\beta)$, and $j\le m$.
Suppose that the probability of getting at least $j$ distinct
values when choosing $m$ i.i.d. random numbers which are uniform
on $\{1,\dots,k\}$ is $1-p(j,k,m)$. Then $Pr(\|P_J^{*m}-U\|>c\beta+\alpha+
p(j,k,m))\le 1/c$ where the probability is over a uniform choice of
$(a_1,\dots,a_k)\in G^k$.
\end{lemma}

{\bf Proof:} Let $I$ be an $m$-tuple $(i_1,\dots,i_m)$ of elements
of $\{1,\dots,k\}$ and $J=(a_1,\dots,a_k)$ be a $k$-tuple of elements
of $G$. Let $\vec\epsilon$ be a vector with $m-j$ elements of
$\{0,1\}$.
Let $S_1$ be the set of $I$ such that $I$ has fewer than $j$ distinct values,
and let $S_2$ be the set of $I$ such that $I$ has at least $j$ distinct
values.

If $I\in S_2$, consider the probability distribution of
$a_{i_1}^{\epsilon_1}a_{i_2}^{\epsilon_2}\dots a_{i_m}^{\epsilon_m}$
where $\epsilon_{\ell}$ is determined by $\vec\epsilon$ if $i_{\ell}=
i_{\ell^{\prime}}$ for some $\ell^{\prime}<\ell$ or $\{i_1,\dots,i_{\ell-1}\}$
has at least $j$ distinct values and where the remaining $j$ values from
$\epsilon_1,\dots,\epsilon_m$ are i.i.d. uniform on $\{0,1\}$. Let
$v(I,J,\vec\epsilon)$ be the variation distance of this probability
distribution from the uniform distribution. By the exercise
\[
\|P_J^{*m}-U\|\le\sum_{I\in S_1}{1\over k^m}1+\sum_{I\in S_2}
\sum_{\vec\epsilon}{1\over k^m}{1\over 2^{m-j}}v(I,J,\epsilon).
\]

Let
\[
G(I,J,\vec\epsilon)=\cases{1&if $v(I,J,\vec\epsilon)\le\alpha$\cr 0&otherwise}
\]
For each $I\in S_2$ and $\vec\epsilon$, the number of $J$ with
$G(I,J,\vec\epsilon)=0$ is no more than $\beta$ times the total number
of $J$ since the family of probability distributions $a_{i_1}^{\epsilon_1}
a_{i_2}^{\epsilon_2}\dots a_{i_m}^{\epsilon_m}$ (where $j$ of the
values $\epsilon_1,\dots,\epsilon_m$ are i.i.d. uniform on $\{0,1\}$
and the rest are determined by $\vec\epsilon$ as previously described)
is $(\alpha,\beta)$-good in variation distance by Lemma~\ref{goodlemma}.
For a given $J$, the number of $I\in S_2$ and $\vec\epsilon$ with
$G(I,J,\vec\epsilon)=0$ may be more than $c\beta$ times the total number
of $I\in S_2$ and $\vec\epsilon$. However, the number of such $k$-tuples
$J$ can be at most $1/c$ times the total number of $k$-tuples $J$.
For the other choices of $J$, we have
\begin{eqnarray*}
\lefteqn{\|P_J^{*m}-U\|}\hspace{0.2in}\\
&\le&\sum_{I\in S_1}{1\over k^m}1+\sum_{(I,\vec\epsilon):
I\in S_2, G(I,J,\vec\epsilon)=1}{1\over k^m}{1\over 2^{m-j}}\alpha+
\sum_{(I,\vec\epsilon): I\in S_2, G(I,J,\vec\epsilon)=0}{1\over k^m}
{1\over 2^{m-j}}1\\
&\le &p(j,k,m)+\alpha+c\beta.
\end{eqnarray*}
The proof of the lemma is complete. \qed

Using Lemma~\ref{puttogether} to prove Theorems~\ref{jtpone} and
\ref{jtptwo} involves finding a bound on $p(j,k,n)$ and choosing $c$
appropriately. The technique to find the bound involves the time
it takes to choose $j$ out of $k$ objects in the coupon collector's
problem. More details may be found in Pak~\cite{pak} and 
in Hildebrand~\cite{mvhjtp}; these sources refer to p. 225 of
Feller~\cite{fellerone}.

By using comparison theorems (Theorem 3 and Proposition 7 of
Hildebrand~\cite{mvhstatprob}), one can extend Theorems~\ref{jtpone}
and \ref{jtptwo} to deal with the cases where
\[
P_J(s)={1-a\over k}|\{i:a_i=s\}|+a\delta_{\{s=e\}}
\]
if $a$ is a constant in the interval $(0,1)$ and $k$ is as in those theorems.
The constant multiple in the expression for $m$ may depend on $a$.

PROBLEM FOR FURTHER STUDY. Can these theorems be extended to the case where
$a=0$ for these values of $k$?

Hildebrand in \cite{mvhjtp} and \cite{mvhstatprob} also considers
some random symmetric lazy random walks and again extends results of
Pak~\cite{pak}. One of these results from \cite{mvhjtp} is the following.
\begin{theorem}
\label{symmetriclazy}
Suppose $X_m=a_{i_1}^{\epsilon_1}\dots a_{i_m}^{\epsilon_m}$ where
$\epsilon_1, \dots, \epsilon_m$ are i.i.d. with $P(\epsilon_i=1)=
P(\epsilon_i=-1)=1/4$ and $P(\epsilon_i=0)=1/2$. Suppose $i_1, \dots, i_m$
are i.i.d. uniform on $\{1,\dots,k\}$ where $k$ is as in Theorem~\ref{jtptwo}.
Given $J\in G^k$, let $Q_{sym}$ be the probability of $X_1$. If
$m=m(n)>(1+\epsilon)(\log_2n)\log(\log_2n)$, then $E(\|Q_{sym}^{*m}-U\|)
\rightarrow 0$ as $n\rightarrow\infty$ where the expectation is over a
uniform choice of $J=(a_1,\dots,a_k)$ from $G^k$.
\end{theorem}

\subsection{Some random random walks on $\grpz_2^d$}

Greenhalgh~\cite{greenhalghpaper} uses some
fairly elementary arguments to examine random random walks on
$\grpz_2^d$, and Wilson~\cite{wilson} uses a binary entropy function
argument to examine these random random walks. The main result of
\cite{wilson} is the following theorem, which we state but do not prove.
\begin{theorem}
\label{wilsontheorem}
Suppose $k>d$.
There exists a function $T(d,k)$ such that the following holds.
Let $\epsilon>0$ be given.
For any choice of $a_1,\dots,a_k$ (each from $\grpz_2^d$), if
$m\le (1-\epsilon)T(d,k)$, then $\|P_{a_1,\dots,a_k}-U\|>1-\epsilon$.
For almost all choices of $a_1,\dots,a_k$, if $m\ge (1+\epsilon)T(d,k)$, then
$\|P_{a_1,\dots,a_k}-U\|<\epsilon$ provided that the Markov chain is
ergodic.
\end{theorem} 

Note that ``for almost all choices'' a property holds
means here that with probability approaching $1$ as 
$d\rightarrow\infty$, the property holds.
Also note that here
$P_{a_1,\dots,a_k}(s)=|\{i:a_i=s\}|/k.$

Some properties of $T(d,k)$ are described in \cite{wilson} and are also
mentioned on p. 321 of Saloff-Coste~\cite{saloffcoste}.

Wilson~\cite{wilson} noted that the upper bound remains valid for
any finite abelian group $G$ provided that $d$ is replaced by $\log_2|G|$.

PROBLEM FOR FURTHER STUDY. Does the expression for the upper bound remain
valid for any finite group $G$ provided that $d$ is replaced by
$\log_2|G|$?

PROBLEM FOR FURTHER STUDY. Relatively little is known about random random
walks on specific families of finite non-abelian groups if 
$k<\log_2|G|$. Indeed, Saloff-Coste~\cite{saloffcoste} (p. 324)
cites a wide-open
problem involving the alternating group with $k=2$.

\end{document}